\documentclass[10pt]{amsart}

\usepackage{amssymb}
\usepackage{stmaryrd}
\usepackage{mathrsfs}
\usepackage{array,subfigure}
\usepackage{graphicx}

\newtheorem{thm}{Theorem}[section]
\newtheorem{prop}[thm]{Proposition}
\newtheorem{cor}[thm]{Corollary}
\newtheorem{lem}[thm]{Lemma}
\newtheorem{defn}[thm]{Definition}
\newtheorem{conj}[thm]{Conjecture}

\newenvironment{xpl}{\refstepcounter{thm} \medskip \noindent {\bf  Example \arabic{section}.\arabic{thm}}}{\hfill$\diamondsuit$\mbox{}\bigskip}

\newcounter{num}
\newenvironment{thmlist}{\begin{list}{(\roman{num})}{\usecounter{num}\setlength{\leftmargin}{25pt}
\setlength{\itemindent}{0pt}\setlength{\labelwidth}{20pt}\setlength{\labelsep}{5pt}\setlength{\itemsep}{0in}}}{\end{list}}

\newcommand{\C}{\mathbb{C}}

\newcommand{\R}{\mathbb{R}}
\newcommand{\Z}{\mathbb{Z}}
\newcommand{\N}{\mathbb{N}}
\newcommand{\Q}{\mathbb{Q}}
\newcommand{\cps}{\mathbb{C}P}

\newcommand{\ol}[1]{\bar{#1}}
\newcommand{\sm}[1]{\scriptscriptstyle{#1}}
\newcommand{\codim}{\operatorname{codim}}
\newcommand{\contr}{\,\lrcorner\,}

\newcommand{\Aut}{\operatorname{Aut}}

\newcommand{\im}{\operatorname{Im}}

\newcommand{\ric}{\operatorname{Ricci}}
\newcommand{\Ric}{\operatorname{Ric}}

\newcommand{\rank}{\operatorname{rank}}

\newcommand{\SF}{\operatorname{SF}}

\newcommand{\Pic}{\operatorname{Pic}}

\newcommand{\Sing}{\operatorname{Sing}}
\newcommand{\inter}{\operatorname{Int}}

\title{Ricci-flat K\"{a}hler metrics on crepant resolutions of K\"{a}hler cones}
\author{Craig van Coevering}
\address{Department of Mathematics, Massachusetts Institute of Technology, 77 Massachusetts Avenue, Cambridge, MA 02139-4307}
\email{craig@math.mit.edu}
\date{May 1, 2008}
\keywords{Calabi-Yau manifold, Sasaki manifold, Einstein metric, Ricci-flat manifold, toric varieties}
\subjclass{Primary 53C25, Secondary 53C55, 14M25 }

\begin{document}

\begin{abstract}
We prove that a crepant resolution $\pi:Y\rightarrow X$ of a Ricci-flat K\"{a}hler cone $X$ admits a complete
Ricci-flat K\"{a}hler metric asymptotic to the cone metric in every K\"{a}hler class in $H^2_c(Y,\R)$.
A K\"{a}hler cone $(X,\ol{g})$ is a metric cone over a Sasaki manifold $(S,g)$, i.e. $X=C(S):=S\times\R_{>0}$
with $\ol{g}=dr^2 +r^2 g$, and $(X,\ol{g})$ is Ricci-flat precisely when $(S,g)$ Einstein of positive scalar
curvature.  This result contains as a subset the existence of ALE Ricci-flat K\"{a}hler metrics on crepant
resolutions $\pi:Y\rightarrow X=\C^n /\Gamma$, with $\Gamma\subset SL(n,\C)$, due to
P. Kronheimer $(n=2)$ and D. Joyce $(n>2)$.

We then consider the case when $X=C(S)$ is toric.  It is a result of A. Futaki, H. Ono, and G. Wang that any
Gorenstein toric K\"{a}hler cone admits a Ricci-flat K\"{a}hler cone metric.
It follows that if a toric K\"{a}hler cone $X=C(S)$ admits a crepant resolution $\pi: Y\rightarrow X$,
then $Y$ admits a $T^n$-invariant Ricci-flat K\"{a}hler metric asymptotic to the cone metric $(X,\ol{g})$ in every
K\"{a}hler class in $H^2_c(Y,\R)$.  A crepant resolution, in this context, is a simplicial fan refining
the convex polyhedral cone defining $X$.  We then list some examples which are easy to construct
using toric geometry.
\end{abstract}

\maketitle

\section{Introduction}

There has been much research recently in constructing examples of Sasaki-Einstein manifolds
(cf.~\cite{BGJ,BG1,MS1,GMSW1,GMSW2,FOW}).
Recall that a Sasaki-Einstein manifold $(S,g)$ is positive scalar curvature Einstein manifold whose metric cone
$(C(S),\ol{g})$, $C(S)=\R_{>0} \times S$ and $\ol{g}=dr^2 +r^2 g$, is Ricci-flat K\"{a}hler.  In all cases besides
$S =S^{2n-1}$, $C(S)=\C^{n}$, the cone has a singularity at the apex.  There has been interest recently in
constructing Ricci-flat K\"{a}hler metrics on resolutions $\pi: Y\rightarrow X$ of the singularity of $X$.
One source of interest in these asymptotically conical Calabi-Yau manifolds is in the AdS/CFT correspondence
(cf.~\cite{MS4,MS3}).  Another motivation is in the construction of new Calabi-Yau manifolds by resolving
conical singularities of a singular Calabi-Yau space (cf~\cite{CanOs,Cha}).

The resolution will necessarily be crepant, and one requires that the metric on $Y$ be asymptotic to the original
Ricci-flat K\"{a}hler cone metric on $X$.  In this article we will give a partial solution to the existence of such
metrics.  Many examples are already known.  In particular,
when $C(S)=\C^{n}/\Gamma$, for $\Gamma\subset SL(n,\C)$ a finite group acting freely on $\C^n \setminus\{o\}$,
such a metric on a resolution of $X$ will be an
ALE Ricci-flat K\"{a}hler metric.  The existence and uniqueness of ALE Ricci-flat K\"{a}hler metrics, in each K\"{a}hler
class, on $X$ has been proved by P. Kronheimer~\cite{Kro} for $n=2$ and by D. Joyce~\cite{Joy1,Joy2} for $n>2$.

Until recently all known examples of Sasaki-Einstein manifolds were quasi-regular, meaning they are orbifold
fibrations over K\"{a}hler Einstein orbifolds.  The 5-dimensional Sasaki-Einstein manifolds $Y^{p,q}$ of
J. Gauntlett, D. Martelli, J. Sparks, and D. Waldram~\cite{GMSW1} provided the first irregular examples,
meaning that they are not simply orbifold fibrations over a K\"{a}hler-Einstein orbifold.
The general existence problem for Sasaki-Einstein metrics on toric Sasaki manifolds has been solved in general in the
beautiful paper of A. Futaki, H. Ono, G. Wang~\cite{FOW}.  In other words, their result implies that any toric
$\Q$-Gorenstein isolated singularity $X$ admits a Ricci-flat K\"{a}hler cone metric.  This will be used as a
source of examples in this article.  Although a crepant resolution does not always exist, it is elementary to
construct examples using toric geometry.

Previous constructions of complete Ricci-flat K\"{a}hler metrics such as those of
G. Tian and S.-T. Yau~\cite{TY1,TY2} constructed metrics asymptotic, in some sense, to a cone over a regular,
or quasi-regular, Sasaki-Einstein manifold.  The present work differs in that the existence of complete Ricci-flat
metrics are proved which are asymptotic to cones over irregular Sasaki-Einstein manifolds.  Some explicit examples of
Ricci-flat K\"{a}hler metrics asymptotic to cones over irregular Sasaki-Einstein manifolds were constructed
in~\cite{MS2}.  The author has also considered the possibility of such metrics on quasi-projective
manifolds in~\cite{vC2}, which is very complementary to this article.

This article considers the following conjecture which first appeared in~\cite{MS3}.
\begin{conj}\label{con:main}
Let $\pi :Y\rightarrow X$ be a crepant resolution of an isolated singularity $X=C(S)$, where $C(S)$ admits a
Ricci-flat K\"{a}hler cone metric.  Then $Y$ admits a unique Ricci-flat K\"{a}hler metric in each K\"{a}hler class
in $H^2(Y,\R)$ that is asymptotic to a cone over the Sasaki-Einstein manifold $(S,g)$.
\end{conj}

We give a partial solution to this conjecture.  We prove the following, where $H_c^2(Y,\R)$ denotes cohomology with
compact supports.
\begin{thm}\label{thm:main}
Let $\pi :Y\rightarrow X$ be a crepant resolution of the isolated singularity of $X=C(S)$, where $C(S)$ admits a
Ricci-flat K\"{a}hler cone metric.  Then $Y$ admits a Ricci-flat K\"{a}hler metric $g$ in each K\"{a}hler class
in $H_c^2(Y,\R)\subset H^2(Y,\R)$ which is asymptotic to the K\"{a}hler cone metric $\ol{g}$ on $X$ as follows.
There is an $R>0$ such that, for any $\delta>0$ and $k\geq 0$,
\begin{equation}\label{eq:asymp}
\nabla^k \left(\pi_* g -\ol{g}\right) =O\left( r^{-2n+\delta-k}\right)\quad\text{on }\{y\in C(S):r(y)>R\},
\end{equation}
where $\nabla$ is the covariant derivative of $\ol{g}$.
\end{thm}
Note that the inclusion of compactly supported cohomology in this case induces an inclusion
$H_c^2(Y,\R)\subset H^2(Y,\R)$.  And $H_c^2(Y,\R)$ is the subset of $H^2(Y,\R)$ whose restriction to
$S\subset X$ vanishes.  Also, if $\omega$ is a K\"{a}hler class in $H_c^2(Y,\R)$, then it has a
d-dimensional neighborhood of K\"{a}hler classes, where $d=\dim H_c^2(Y,\R)$.  Thus the theorem gives families of
Ricci-flat metrics.  Note also that $d$ is the number of prime divisors in the exceptional set $E=\pi^{-1}(o)$.

It is also useful to consider partial crepant resolutions $\pi:Y\rightarrow X$ where $Y$ has only orbifold
singularities.  The proof of Theorem~\ref{thm:main} is valid without modification in this case also.  Many
of the examples of Sasaki-Einstein manifolds $S$ have associated Ricci-flat K\"{a}hler cones $X=C(S)$ which do
not admit crepant resolutions, but nonetheless admit such a partial crepant resolution.  This is true of
some of the examples constructed via hypersurface singularities in~\cite{BG1} and~\cite{BGJ}, while some examples
do admit crepant resolutions.

Theorem~\ref{thm:main} solves a large portion of Conjecture~\ref{con:main}.  But it is instructive to consider
what it excludes.  If $\pi :Y\rightarrow X$ is a small resolution, i.e. $\codim_{\C}(E)>1$, where $E=\pi^{-1}(o)$
is the exceptional set, then there are no K\"{a}hler classes in $H_c^2(Y,\R)$.  In particular, consider
the conifold $X=\{z_0^2 +z_1^2 +z_2^2 +z_3^2 =0\}\subset\C^4$ which is the cone over $S^2 \times S^3$.
It has the structure of a Ricci-flat K\"{a}hler cone if $S^2 \times S^3$ is given the homogeneous
Sasaki-Einstein metric.  Then $X$ admits a crepant resolution $\pi:Y\rightarrow X$, where
$Y$ is the total space of $\mathcal{O}(-1)\oplus\mathcal{O}(-1)\rightarrow\cps^1$.  The exceptional set
is $\cps^1 =\pi^{-1}(o)$.  Nevertheless, it is well known that $Y$ admits a complete Ricci-flat K\"{a}hler
metric (cf.~\cite{CanOs}).

After proving Theorem~\ref{thm:main} we will consider the toric case in more detail.  In this case $X=C(S)$
is a Gorenstein toric K\"{a}hler cone which admits a toric Ricci-flat K\"{a}hler cone metric by the
results of~\cite{FOW}.  In this case a crepant resolution $\pi:Y\rightarrow X$ is toric, and $Y$ is described
explicitly by a nonsingular simplicial fan $\tilde{\Delta}$ refining the convex polyhedral cone $\Delta$ defining $X$.
A K\"{a}hler class in $H_c^2(Y,\R)$ is characterized by a \emph{compact} strictly convex
support function on $\tilde{\Delta}$.  This is a strictly convex support function on $\tilde{\Delta}$ which vanishes on the
rays defining $\Delta$.  We prove the following.
\begin{cor}\label{cor:main}
Let $\pi:Y\rightarrow X$ be a crepant resolution of a Gorenstein toric K\"{a}hler cone $X$.
Suppose the fan $\tilde{\Delta}$ defining $Y$ admits a compact strictly convex support function.
Then $Y$ admits a Ricci-flat K\"{a}hler metric $g$ which is asymptotic to $(C(S),\ol{g})$ as in
(\ref{eq:asymp}).  Furthermore, $g$ is invariant under the compact $n$-torus $T^n$.
\end{cor}
As above, if a Ricci-flat K\"{a}hler metric exists, then there is a d-dimensional, $d=\dim H_c^2(Y,\R)$, family
of such metrics.  Here $d$ is the number of lattice points in the interior of a polytope $P_{\Delta}$ which is the
intersection of the cone $\Delta$ defining $X$ with a hyperplane.  Thus, although crepant resolutions are generally
not unique, $d$ is invariant.  A crepant resolution of $X$ is characterized by a basic lattice triangulation
of $P_{\Delta}$.  When $n=3$ such a triangulation always exists.

In the final section we give some examples.  These are easily described by the toric geometry of the resolution
$\pi:Y\rightarrow X$ in the toric case.  Many more examples are constructed in~\cite{vC3} using toric geometry
and by resolving hypersurface singularities.   Recently the author has come up with a proof which removes the $\delta>0$
from the convergence in (\ref{eq:asymp}) and thus gives the sharp convergence.  This will appear in a subsequent article.

\section{Sasaki manifolds}

We review some of the properties of Sasaki manifolds.  We are primarily interested in K\"{a}hler cones.  But a K\"{a}hler cone
is a cone over a Sasaki manifold, and much research has been done recently on Sasaki-Einstein manifolds (cf.~\cite{BG,BG2}).
\begin{defn}
 A Riemannian manifold $(S,g)$ of dimension $2n-1$ is Sasakian if the metric cone $(C(S),\ol{g})$, $C(S)=\R_{>0}\times S$
and $\ol{g}=dr^2 +r^2 g$, is K\"{a}hler.
\end{defn}

Set $\tilde{\xi}=J(r\frac{\partial}{\partial r})$, then $\tilde{\xi}-iJ\tilde{\xi}$ is a holomorphic vector field
on $C(S)$.  The restriction $\xi$ of $\tilde{\xi}$ to $S=\{r=1\}\subset C(S)$
is the \emph{Reeb vector field} of $S$, which is a Killing vector field.
If the orbits of $\xi$ close, then it defines a locally free $U(1)$-action on $S$.  If the $U(1)$-action is free,
then the Sasaki structure is said to be \emph{regular}.  If there are non-trivial stablizers then the Sasaki
structure is \emph{quasi-regular}.  If the orbits do not close the Sasaki structure is \emph{irregular}.

Let $\eta$ be the dual 1-form to $\xi$ with respect to $g$.  Then
\begin{equation}\label{eq:cont}
\eta = (2d^c \log r)|_{r=1},
\end{equation}
where $d^c=\frac{i}{2}(\ol{\partial}-\partial)$.  Let $D=\ker\eta$.  Then $d\eta$ in non-degenerate on $D$
and $\eta$ is a contact form on $S$.  Furthermore, we have
\begin{equation}
d\eta(X,Y) =2g(\Phi X,Y), \quad\text{for }X,Y\in D_x, x\in S,
\end{equation}
where $\Phi|_{D_x}$ is the restriction of the complex structure $J$ on $C(S)$, to $D_x$, and $\Phi(\xi)=0$.  Thus
$(D,J)$ is a strictly pseudo-convex CR structure on $S$.
We will denote the Sasaki structure on $S$ by $(g,\xi,\eta,\Phi)$.
It follows from (\ref{eq:cont}) that the K\"{a}hler form of $(C(S),\ol{g})$ is
\begin{equation}\label{eq:kahler-form}
 \omega=\frac{1}{2}d(r^2 \eta)=\frac{1}{2}dd^c r^2.
\end{equation}
Thus $\frac{1}{2}r^2$ is a K\"{a}hler potential for $\omega$.

There is a 1-dimensional foliation $\mathscr{F}_\xi$ generated by the Reeb vector field $\xi$.  Since the leaf
space is identical with that generated by $\tilde{\xi}-iJ\tilde{\xi}$ on $C(S)$, $\mathscr{F}_\xi$ has
a natural transverse holomorphic structure.  And $\omega^T =\frac{1}{2}d\eta$ defines a K\"{a}hler form on the
leaf space.  We denote the transverse K\"{a}hler metric by $g^T$.
Note that when the Sasaki structure on $S$ is regular (resp. quasi-regular), the leaf space of
$\mathscr{F}_\xi$ is a K\"{a}hler manifold (resp. orbifold).

A p-form $\alpha\in\Omega^p(S)$ on $S$ is said to be \emph{basic} if
\begin{equation}
\xi\contr\alpha =0\quad\text{and}\quad \mathcal{L}_\xi \alpha =0.
\end{equation}
The basic p-forms are denoted by $\Omega_B^p(S)$, where the foliation $\mathscr{F}_\xi$ on $S$ must be fixed.
One easily checks that $\Omega_B^*$ is closed under the exterior derivative.  So there is a transversal de Rham complex
which can be used to calculate the basic cohomology $H^*_B(S)$.

The foliation $\mathscr{F}_\xi$ associated to a Sasaki structure has a transverse holomorphic structure, so there is a
splitting $\Omega_B ^k =\bigoplus_{p+q=k} \Omega_B^{p,q}$ of complex forms into types.  And the exterior derivative on
basic forms splits into $d=\partial +\ol{\partial}$, where $\partial$ has degree $(1,0)$ and $\ol{\partial}$ has degree $(0,1)$.
Thus we have as well the basic Dolbeault complex and the basic Dolbeault cohomology groups $H^{p,q}_B(S)$.

Furthermore, the foliation has a transverse K\"{a}hler structure, and the usual Hodge theory for K\"{a}hler manifolds carries over.
In particular, we have the Hodge decomposition $H^k_B(S,\C)=\bigoplus_{p+q=k} H^{p,q}_B(S)$ and the representation of basic cohomology classes by harmonic forms.  It is also useful to know that the $\partial\ol{\partial}$-lemma holds for basic forms as it does on K\"{a}hler
manifolds.  Thus if $\phi\in\Omega^{1,1}_B$ is exact, then there is a basic $f\in C_B^{\infty}$ with $\phi=i\partial\ol{\partial}f$ and
$f$ can be taken to be real if $\phi$ is.  See the monograph~\cite{BG2} for a survey of these results.

We will consider deformations of the transverse K\"{a}hler structure.  Let $\phi\in C^\infty_B(S)$ be a smooth basic
function.  Then set
\begin{equation}
\tilde{\eta} =\eta +2d^c_B \phi.
\end{equation}
Then
\[ d\tilde{\eta} =d\eta +2d_B d^c_B \phi =d\eta +2i\partial_B \ol{\partial}_B \phi. \]
For sufficiently small $\phi$, $\tilde{\eta}$ is a non-degenerate contact form in that $\tilde{\eta}\wedge d\tilde{\eta}^n$
is nowhere zero.  Then we have a new Sasaki structure on $S$ with the same Reeb vector field $\xi$, transverse
holomorphic structure on $\mathscr{F}_\xi$, and holomorphic structure on $C(S)$.  This Sasaki structure has
transverse K\"{a}hler form $\tilde{\omega}^T=\omega^T +i\partial_B \ol{\partial}_B \phi$.  One can show~\cite{FOW}
that if
\[\tilde{r} =r\exp{\phi},\]
then $\tilde{\omega}=\frac{1}{2}dd^c \tilde{r}^2$ is the K\"{a}hler form on $C(S)$ associated to the
transversally deformed Sasaki structure.

\begin{prop}\label{prop:ricci}
Let $(S,g)$ be a $2n-1$-dimensional Sasaki manifold.  Then the following are equivalent.
\begin{thmlist}
 \item $(S,g)$ is Sasaki-Einstein with the Einstein constant being necessarily $2n-2$.

 \item $(C(S),\ol{g})$ is a Ricci-flat K\"{a}hler.

 \item The K\"{a}hler structure on the leaf space of $\mathscr{F}_\xi$ is K\"{a}hler-Einstein with Einstein constant $2n$.
\end{thmlist}
\end{prop}
This follows from elementary computations.  In particular, the equivalence of (i) and (iii) follows from
\begin{equation}\label{eq:ricci}
 \Ric_g(\tilde{X},\tilde{Y})=(\Ric^T -2g^T)(X,Y),
\end{equation}
where $\tilde{X},\tilde{Y}\in D$ are lifts of $X,Y$ in the local leaf space; $g^T$ and $\Ric^T$ are the metric and Ricci
tensor of the transversal K\"{a}hler structure.

Given a Sasaki structure we can perform a $D$-homothetic transformation to get a new Sasaki structure.  For $a>0$ set
\begin{gather}\label{eq:d-homo}
 \eta'=a\eta,\quad \xi'=\frac{1}{a}\xi,\\
 g'=ag^T +a^2\eta\otimes\eta =ag+(a^2-a)\eta\otimes\eta.\\
\end{gather}
Then $(g',\xi',\eta',\Phi)$ is a Sasaki structure with the same holomorphic structure on $C(S)$, and with
$r'=r^a$.

\begin{prop}\label{prop:CY-cond}
The following necessary conditions for $S$ to admit a deformation of the transverse K\"{a}hler structure to a Sasaki-Einstein
metric are equivalent.
\begin{thmlist}
 \item $c_1^B =a[d\eta]$ for some positive constant $a$.

 \item $c_1^B >0$, i.e. represented by a positive $(1,1)$-form, and $c_1(D)=0$.

 \item For some positive integer $\ell>0$, the $\ell$-th power of the canonical line bundle
 $\mathbf{K}^{\ell}_{C(S)}$ admits a nowhere vanishing section $\Omega$ with
 $\mathcal{L}_\xi \Omega =i n\Omega$.
\end{thmlist}
\end{prop}
\begin{proof}
 Let $\rho$ denote the Ricci form of $(C(S),\ol{g})$ and $\rho^T$ the Ricci form of $\Ric^T$, then easy computation shows that
 \begin{equation}\label{eq:ricci-cone}
  \rho =\rho^T - 2n\frac{1}{2}d\eta.
 \end{equation}
If (i) is satisfied, there is a $D$-homothety so that $[\rho^T]= 2n[\frac{1}{2}d\eta]$ as basic classes.
Thus there exists a smooth function $h$ with $\xi h=0=r\frac{\partial}{\partial r}h$ and
\begin{equation}
 \rho =i\partial\ol{\partial}h.
\end{equation}
This implies that $e^h \frac{\omega^{n}}{n!}$, where $\omega$ is the K\"{a}hler form of $\ol{g}$, defines
a flat metric $|\cdot|$ on $\mathbf{K}_{C(S) }$.  Parallel translation defines a multi-valued section which defines
a holomorphic section $\Omega$ of $\mathbf{K}^{\ell}_{C(S)}$ for some integer $\ell>0$ with
$|\Omega |=1$.  Then we have
\begin{equation}\label{eq:hol-form}
 \left(\frac{i}{2}\right)^{n}(-1)^{\frac{n(n-1)}{2}}\Omega\wedge\ol{\Omega} =e^h\frac{1}{n!}\omega^{n}.
\end{equation}
From the invariance of $h$ and the fact that $\omega$ is homogeneous of degree 2, we see that
$\mathcal{L}_{r\frac{\partial}{\partial r}}\Omega= n\Omega$.

Conversely, if (iii) holds, then we have (\ref{eq:hol-form}) for some $h\in C^\infty(C(S))$.  Then since
$\omega$ is homogeneous of degree 2 and $\mathcal{L}_{r\frac{\partial}{\partial r}}\Omega= n\Omega$, it follows that
$\xi h=0=r\frac{\partial}{\partial r}h$.  And the above arguments show that $c_1^B =\frac{1}{2\pi}[\rho^T]=\frac{n}{\pi}[d\eta]$.

The equivalence of (i) and (ii) is easy (cf.~\cite{FOW} Proposition 4.3).
\end{proof}

\begin{xpl}\label{xpl:Sasak-st}
Let $Z$ be a complex manifold (or orbifold) with a negative holomorphic line bundle (respectively V-bundle) $\mathbf{L}$.
If the total space of $\mathbf{L}^\times$, $\mathbf{L}$ minus the zero section, is smooth, then the $U(1)$-subbundle
$S\subset\mathbf{L}^\times$ has a natural regular (respectively quasi-regular) Sasaki structure.
Let $h$ be an Hermitian metric on $\mathbf{L}$ with negative curvature.  If in local holomorphic coordinates we define
$r^2 =h|z|^2$, where $z$ is the fiber coordinate, then $\omega =\frac{1}{2}dd^c r^2$ is the K\"{a}hler form on
$\mathbf{L}^\times$ of a K\"{a}hler cone metric.  And $S=\{z\in\mathbf{L}^\times : r(z)=1\}$ has the induced Sasaki structure.
Conversely, it can be shown that every regular (respectively quasi-regular) Sasaki structure arises from this construction
(cf.~\cite{BG}).
\end{xpl}

\section{Crepant resolutions}

Let $C(S)$ be a K\"{a}hler cone.  Note that \textit{a priori} $C(S)$ does not contain the vertex, but
$X=C(S)\cup\{o\}$ can be made into a complex space in a unique way.  The Reeb vector field $\xi$ generates
a 1-parameter subgroup of the automorphism group $\Aut(S)$ of the Sasaki manifold $S$.  Since $\Aut(S)$ is compact,
the closure of this subgroup is a torus $T^k \subset\Aut(S)$.  Here $\rank(S):= k$.  Choose a vector field $\zeta$ in the
integral lattice of the Lie algebra $\mathfrak{t}$ of $T^k$, $\zeta\in\Z_T \subset\mathfrak{t}$, and such that $\eta(\zeta)>0$ on $S$.
Then it is not difficult to show that there is a quasi-regular Sasaki structure $(\tilde{g},\zeta,\tilde{\eta},\tilde{\Phi})$
with the same CR-structure $(D,J)$ and with Reeb vector field $\zeta$ (cf.~\cite{BGS2}).
And the $U(1)$-action on $S$ generated by $\zeta$ extends to
an holomorphic $\C^*$-action on $C(S)$.  Then the quotient $C(S)/\C^* = S/ U(1)$ is a K\"{a}hler orbifold $Z$, and
$C(S)$ is the total space, minus the zero section, of an orbifold bundle $\iota:\mathbf{L}\rightarrow Z$
(cf.~\cite{BG}).  The bundle $\mathbf{L}$ is negative.  There is a metric $h$ on $\mathbf{L}$, so that
$\tilde{r}^2 =h|z|^2$ locally, where $z$ is the fiber coordinate.  And the K\"{a}hler form on $C(S)$ for the
Sasaki structure $(\tilde{g},\zeta,\tilde{\eta},\tilde{\Phi})$ is $\frac{1}{2}i\partial\ol{\partial}\tilde{r}^2$ as in
(\ref{eq:kahler-form}).

Let $W$ be the total space of $\mathbf{L}$.  Then $\tilde{r}^2$ is strictly
plurisubharmonic away from $Z\subset W$, and hence $W$ is a 1-convex space.  In other words, $W$ is exhausted by
strictly pseudo-convex domains $\{\tilde{r}^2 <c \}\subset W$, for $c>0$.  Then as in~\cite{Gra} $W$ is holomorphically
convex, and we have the Remmert reduction of $W$.  That is, there exists a Stein space $X$ and an holomorphic map
$\sigma: W\rightarrow X$, which contracts the maximal compact analytic set $Z\subset W$ and is a biholomorphism
outside $Z$.  Thus $X=C(S)\cup\{o\}$ is a complex space. Furthermore $X$ is normal, and the Riemann extension theorem
shows $i_* \mathcal{O}_{C(S)} =\mathcal{O}_X$, where $i: C(S)\rightarrow X$ is the inclusion.  Thus
$X$ is independent of the above choices.

Note that $X=C(S)\cup\{o\}$ is a Stein space.  And if $\pi:Y\rightarrow X$ is any resolution of $o\in X$, then
$Y$ is 1-convex.  It is actually known~\cite{OrVer} that $X=C(S)\cup\{o\}$ is an affine variety.  See also~\cite{vC3}
for a succinct proof.

Recall that a singularity $x\in X$ is \emph{rational} if $(R^i \pi_* \mathcal{O}_Y)_x =0$, for $i>0$, where $\pi:Y\rightarrow X$
is a resolution of singularities.  One can show that this is independent of the resolution.

Suppose $o\in X$ is an isolated singularity.  Then we have a simple criterion for
rationality (cf.~\cite{Bur} and~\cite{Lau}).
\begin{prop}\label{prop:ration}
Let $\Omega$ be a holomorphic $n$-form defined, and nowhere vanishing, on a deleted neighborhood of $o\in X$.
Then $o\in X$ is rational if and only if
\begin{equation}\label{eq:ration}
\int_U \Omega\wedge\ol{\Omega}<\infty,
\end{equation}
for $U$ a sufficiently small neighborhood of $o\in X$.
\end{prop}
Note that if (\ref{eq:ration}) is satisfied for $\Omega$, then it is satisfied for all holomorphic $n$-forms defined in a neighborhood
of $o\in X$.  And for any such form $\pi^*\Omega$ extends to a holomorphic form on $\tilde{U}=\pi^{-1}(U)$.

Let $\omega_X$ denote the dualizing sheaf of $X$.  Then we have $\omega_X \cong i_*(\mathcal{O}(\mathbf{K}_{C(S)}))$,
where $i:C(S)\rightarrow X$ is the inclusion, as the codimension on $\Sing(X)=\{o\}\subset X$ is greater than 2.
Recall that $X$ is said to be p-Gorenstein if $\omega_X^{[p]}:=i_*(\omega_{C(S)}^{\otimes p})$ is locally free for
$p\in\N$, and $X$ is $\Q$-Gorenstein if it is p-Gorenstein for some p.  We will call $X$ Gorenstein if it is
1-Gorenstein.

Suppose $X$ is $\Q$-Gorenstein.  A resolution $\pi:Y\rightarrow X$ is said to be \emph{crepant} if
\begin{equation}
\pi^* \omega_X =\omega_Y =\mathcal{O}(\mathbf{K}_Y).
\end{equation}

\begin{prop}\label{prop:resol}
Let $X=C(S)$ be the K\"{a}hler cone of a Sasaki manifold $S$ satisfying Proposition~\ref{prop:CY-cond},
e.g. $S$ is Sasaki-Einstein.  Then $X$ is $\Q$-Gorenstein, and $o\in X$ is a rational singularity.

Suppose $X$ admits a crepant resolution $\pi:Y\rightarrow X$.
If $H_1(Y,\Z)=0$, which is always the case in dimension 3, then $X$ is Gorenstein.
\end{prop}
\begin{proof}
By Proposition~\ref{prop:CY-cond} There exits a section $\Omega_p \in\Gamma(\mathbf{K}_{C(S)}^{\otimes p})$.
The Riemann extension theorem shows
that $\omega_X^{[p]} =i_*(\mathcal{O}(\mathbf{K}_{C(S)}^{\otimes p}))$ is locally free, and in fact trivial.
Thus $X$ is $\Q$-Gorenstein.

Note that the conditions of Proposition~\ref{prop:CY-cond} imply that $\pi_1(S)$ is finite.  Indeed, the
transversal Ricci form $\ric^T \in [a\omega^T] ,\ a>0,$ where $a\omega^T$ is a positive basic $(1,1)$ class.
By the transverse version of the Calabi-Yau theorem there is a transversal K\"{a}hler deformation to a Sasaki
structure with $\Ric^T >0$.  Then after a possible $D$-homothetic transformation,  equation (\ref{eq:ricci})
shows that one can obtain a Sasaki metric with $\Ric_g >0$.  Then the claim follows by Meyer's theorem.

The universal cover $\ol{S}$ of $S$ is finite, and we have a
finite unramified morphism $g :\ol{X}\rightarrow X$, where $\ol{X}=C(\ol{S})\cup\{o\}$.
The holomorphic form $\Omega$ on $C(\ol{S})$ from Proposition~\ref{prop:CY-cond} is easily seen to
satisfy (\ref{eq:ration}).  In fact, the proof of Proposition~\ref{prop:ration} shows that $\Omega$ extends
to a regular form on any resolution of $\ol{X}$.  It is well known that the image of a finite morphism $X$ must also have
rational singularities~\cite[Prop. 5.13]{KolMor}.

By assumption $\pi^*\Omega_r$ is a nonvanishing section of
$\mathbf{K}_Y^{\otimes r}$.  One can prove using the definition of a rational singularity that $\Pic Y =H^2(Y,\Z)$ (cf.~\cite{vC3}),
which is free by assumption.  Thus $\mathbf{K}_Y$ is trivial and has a nowhere vanishing section $\Omega$, and its restriction to
$i_*(\mathcal{O}(\mathbf{K}_{C(S)}))$ defines a nonvanishing section of $\omega_X$.

It is a result of N. Shepherd-Barron than a crepant resolution of an isolated canonical 3-fold singularity is in fact
simply connected.
\end{proof}

\section{Approximate metric}

Let $X=C(S)\cup\{o\}$ be a K\"{a}hler cone.  Suppose $\pi:Y\rightarrow X$ is a resolution of $o\in X$.
We will denote the pull back $\pi^* r$ of the radius function $r$ on $C(S)$ to $Y$ by $r$ also.
Let $\overline{Y}=\{y\in Y :r(y)\leq 1\}\subset Y$.  Then $H^2_c(Y,\R)\cong H^2(\overline{Y},S,\R)$ and the cohomology
sequence gives
\begin{equation}
\cdots\rightarrow H^1(S,\R)\rightarrow H^2_c(Y,\R)\rightarrow H^2(Y,\R)\rightarrow H^2(S,\R)\rightarrow\cdots.
\end{equation}
Suppose that $S$ satisfies Proposition~\ref{prop:CY-cond}, then $H^1(S,\R)=\{0\}$ by the argument in Proposition~\ref{prop:resol}.
Thus we have an inclusion $H^2_c(Y,\R)\subset H^2(Y,\R)$.  In fact, one can prove with some more work that
$0\rightarrow H^2_c(Y,\R)\rightarrow H^2(Y,\R)\rightarrow H^2(S,\R)\rightarrow 0$ is exact (cf.~\cite{vC3}).

We prove that the restriction in Theorem~\ref{thm:main} to K\"{a}hler classes in $H_c^2(Y,\R)$ is in some
sense necessary.
\begin{prop}\label{prop:asymp-kah}
Let $\pi:Y\rightarrow X$ be a resolution of the K\"{a}hler cone $X=C(S)$.  Let $g$ be a K\"{a}hler metric on
$Y$ with K\"{a}hler form $\omega$.  Suppose
\begin{equation}\label{eq:asymp-kah}
\|\pi_* g -\ol{g}\|_{\ol{g}} =O\left(r^{-\alpha}\right),
\end{equation}
where $\ol{g}$ is the cone metric on $C(S)$.  If $\alpha >2$, then $[\omega]\in H_c^2(Y,\R)$.
\end{prop}
\begin{proof}
Let $\ol{\omega}=\frac{1}{2}dd^c r^2$, and set $\beta =\omega -\ol{\omega}$.  Let
$i_a : S\subset C(S)$, for $a>0$, be the inclusion as the set $\{r=a\}\subset C(S)$.
In the following $\gamma\in\Omega^2(S)$ is an arbitrary 2-form, $g_1$ is the Sasaki metric on $S$, and
$g_r =a^2 g_1$ is the metric on $S$ induced by $i_a$:
\begin{equation}\label{eq:asymp-int}
\begin{split}
\int_S i_a^*\beta\wedge\star_{\sm g_1}\gamma & = \int_S \langle i_a^*\beta,\gamma\rangle\mu_{\sm g_1} \\
                                             & \leq \int_S \|i_a^*\beta\|_{\sm g_1}\|\gamma\|_{\sm g_1}\mu_{\sm g_1} \\
                                             & = \int_S a^2 \|i_a^*\beta\|_{\sm g_r}\|\gamma\|_{\sm g_1}\mu_{\sm g_1}.
\end{split}
\end{equation}
And $\|i_a^*\beta\|_{\sm g_r}\leq i_a^* \|\beta\|_{\sm{\ol{g}}}$.  By (\ref{eq:asymp-kah}) there is a constant $C>0$
so that
\begin{equation}
\int_S a^2 \|i_a^*\beta\|_{\sm g_r}\|\gamma\|_{\sm g_1}\mu_{\sm g_1}\leq C\int_S a^{-\alpha+2}\|\gamma\|_{\sm g_1}\mu_{\sm g_1}
\rightarrow 0,\text{  as }a\rightarrow 0.
\end{equation}
If $\star_{\sm g_1}\gamma$ is closed, then the integral on the left of (\ref{eq:asymp-int}) is independent of $a>0$.
\end{proof}

This has consequences in the case of small resolutions.
\begin{cor}
Suppose $\pi:Y\rightarrow X$ is a small resolution.  And $g$ is an asymptotically conical metric on $Y$, meaning
that $g$ satisfies (\ref{eq:asymp-kah}) for some $\alpha>0$.  Then $\alpha\leq 2$.
\end{cor}
\begin{proof}
Suppose $\alpha>2$.  Then by Proposition~\ref{prop:asymp-kah} the K\"{a}hler form $\omega$ satisfies
$[\omega]\in H^2_c(Y,\R)$.  Thus $[\omega]$ is Poincar\'{e} dual an element of $H_{2n-2}(Y,\R)$.  But
$Y$ is homotopically equivalent to $E=\pi^{-1}(o)$, and $\dim_{\C} E <n-1$.  So $H_{2n-2}(Y,\R)=\{0\}$.
\end{proof}

Indeed, there is an asymptotically conical Calabi-Yau metric on the small resolution
$\pi:Y\rightarrow X$, where $Y$ is the total space of $\mathcal{O}(-1)\oplus\mathcal{O}(-1)\rightarrow\cps^1$
and $X=\{z_0^2 +z_1^2 +z_2^2 +z_3^2 =0\}\subset\C^4$, constructed in~\cite{CanOs}.  And for this metric one has
$\alpha=2$.

Suppose that $C(S)$ is a K\"{a}hler cone which satisfies Proposition~\ref{prop:CY-cond}.
And suppose $\pi:Y\rightarrow X$ is a resolution.
\begin{lem}\label{lem:app-form}
Let $\omega'$ be a K\"{a}hler metric on $Y$ whose cohomology class $[\omega']\in H^2_c(Y,\R)$.  Then
there exists a K\"{a}hler metric $\omega_0$ on $Y$ with $[\omega_0]=[\omega']$ and for some $r_0 >0$ on
$Y_{r_0} =\{y\in Y: r(y)\geq r_0\}$ $\omega_0$ restricts to $\pi^*\omega$, the pull-back of the K\"{a}hler cone metric.
\end{lem}
\begin{proof}
Let $\{D_i\}$ be the prime divisors in the exceptional set $E=\pi^{-1}(o)$.  Since $H_{2n-2}(Y,\R)$ is generated
by the fundamental classes of the $D_i$, $[\omega']$ is Poincar\'{e} dual to $\sum_i a_i D_i$ for $a_i \in\R$.
Thus there exists a compactly supported closed $(1,1)$-form $\beta$ Poincar\'{e} dual to $\sum_i a_i D_i$, so
$[\beta]=[\omega']$.

Proposition~\ref{prop:resol} implies that $o\in X$ is a rational singularity, so $(R^i \pi_* \mathcal{O}_Y)_o =0$ for $i>0$.
And because $X$ is a Stein space
$H^j(X,R^i \pi_* \mathcal{O}_Y)=0$, for $j>0$.  Thus the Leray spectral sequence implies that $H^p(Y,\mathcal{O}_Y)=0$ for
$p>0$.

Let $\alpha\in\mathcal{A}^1(Y)$ be a smooth 1-form with $d\alpha =\omega'-\beta$.  Then if $\alpha =\alpha^{1,0} +\alpha^{0,1}$
is the decomposition into types $\alpha^{0,1} =\overline{\alpha^{1,0}}$, and $\ol{\partial}\alpha^{0,1} =0$.  So there exists a smooth,
complex valued, function $\gamma$ with $\ol{\partial}\gamma =\alpha^{0,1}$, since $H^1_{\ol{\partial}}(Y)=H^1(Y,\mathcal{O}_Y)=0$.
And it easily follows that $\mathbf{i}\partial\ol{\partial}(2\im{\gamma})=\omega'-\beta$.

Denote by $f=\frac{r^2}{2}$ the K\"{a}hler potential of the Ricci-flat K\"{a}hler cone metric on $X$.
Choose $0< a_1 <a_2$.  And let $\nu:\R_{>0}\rightarrow\R$ be a smooth function with $\nu(x)=x$ for $x>\frac{a_2^2}{2}$,
$\nu'(x),\nu''(x)\geq 0$ for $\frac{a_1^2}{2}<x<\frac{a_2^2}{2}$, and $\nu(x)=c$, a constant, for $x<\frac{a_1^2}{2}$.  Then
$i\partial\ol{\partial}(\nu\circ f)\geq 0$ and extends to a form on $Y$.  Now choose $b_1,b_2$ with $\frac{a_2^2}{2}<b_1 <b_2$.
Let $\phi$ be a non-negative function of $r$ with $\phi(r)=1$ for $r<b_1$ and $\phi(r)=0$ for $r>b_2$.  Then define
$u=2\phi\im{\gamma}$ and
\begin{equation}
 \omega_0 =\beta +i\partial\ol{\partial}u + Ci\partial\ol{\partial}(\nu\circ f), \text{  for }C>0.
\end{equation}
For $C>0$ sufficiently large this gives the required metric.
\end{proof}

\section{Monge-amp\`{e}re equation}\label{sect:monge}

In this section we prove the existence of the complete Ricci-flat metric, and its asymptotic properties, in
Theorem~\ref{thm:main}.  Basically the arguments in~\cite{TY1} and ~\cite{TY2} for the case of Ricci-flat metrics on
quasi-projective manifolds work in this situation, but in this situation we are able to fix the asymptotics of the
metric more precisely.

Suppose $(C(S),\omega)$ is a Ricci-flat K\"{a}hler cone, and $\pi:Y\rightarrow X$ is a crepant resolution.
There is a holomorphic n-form $\Omega$ on $C(S)$ satisfying (\ref{eq:hol-form}) with $h$ constant.
Thus there is a $c\in\C$ so that
\begin{equation}
c\Omega\wedge\ol{\Omega} =\omega^{n}.
\end{equation}
Let $\Omega$ also denote the extension of $\pi^* \Omega$ to a nowhere vanishing n-form on $Y$.

Define a real valued function
\begin{equation}\label{eq:ricci-pot}
 f=\log\left(\frac{c\Omega\wedge\ol{\Omega}}{\omega_0^n}\right),
\end{equation}
where $\omega_0$ is the K\"{a}hler form of Lemma~\ref{lem:app-form}.
Then $i\partial\ol{\partial}f =\ric(\omega_0)$, and after possibly adding a constant to $f$,
$f$ vanishes on $Y_{b_2} =\{y\in Y : r(y)> b_2\}$.  The existence of of a Ricci-flat
K\"{a}hler metric on $Y$ is equivalent to a solution to the following Monge-Amp\`{e}re equation:
\begin{equation}\label{eq:monge-amp}
\begin{cases}
 \left(\omega_0 +i\partial\ol{\partial}\phi\right)^n =e^f \omega_0^n,\\
 \omega_0 +i\partial\ol{\partial}\phi >0.
\end{cases}
\end{equation}

For the proof of the following see~\cite{TY2}, Proposition 4.1.  Note that the proof makes use of the boundedness of the curvature tensor
$\|R(g_0)\|<\infty$, where $g_0$ is the metric associated to $\omega_0$, and of the covariant derivative of the scalar curvature
$\|ds_{g_0}\|<\infty$.  The proof also makes use of some analysis developed in~\cite{CheY}.
\begin{prop}\label{prop:monge}
Let $\omega_0$ be the K\"{a}hler form defined in Lemma~\ref{lem:app-form}.  Then there is a unique solution $\phi$ to (\ref{eq:monge-amp})
such that $\phi(y)$ converges uniformly to zero as $y$ goes to infinity, and there is a constant $c>1$ so that
$c^{-1}\omega_0 <\omega_0 +i\partial\ol{\partial}\phi <c\omega_0$.  It follows that $\tilde{\omega}=\omega_0 +i\partial\ol{\partial}\phi$ is a complete
Ricci-flat K\"{a}hler metric on $Y$.
\end{prop}

We now prove that the metric $\tilde{g}$ with K\"{a}hler form $\tilde{\omega}$ of Proposition~\ref{prop:monge} is
asymptotic to the Ricci-flat K\"{a}hler cone metric $(C(S),\ol{g})$ as stated in Theorem~\ref{thm:main}.

\begin{lem}\label{lem:conv1}
Let $\phi$ be the solution to (\ref{eq:monge-amp}) given in Proposition~\ref{prop:monge}.
 For any $\delta>0$ there are constants $C,C_{\delta}>0$ so that
\begin{equation}\label{eq:conv}
 -C_\delta (1+r^2(y))^{-n+1}(\log r(y))^\delta \leq\phi(y)\leq C(1+r^2(y))^{-n+1},\quad\text{for  } y\in Y_{r_0},
\end{equation}
where $Y_{r_0}$ is as in Lemma~\ref{lem:app-form}.
\end{lem}
\begin{proof}
 For $r>r_0$ we have $(\omega +dd^c\phi)^n =\omega^n$, where $\omega=\frac{1}{2}dd^c(r^2)=rdr\wedge\eta +\frac{1}{2}r^2 d\eta$.
Set $\rho=Kr^{-2n+2}$.  Then
\begin{equation}
 dd^c \rho  =2(n-1)^2 Kr^{-2n} rdr\wedge\eta -(n-1)Kr^{-2n+2}d\eta,
\end{equation}
so
\begin{equation}
 \omega +dd^c\rho =(1+2(n-1)^2 Kr^{-2n})rdr\wedge\eta +(1-2(n-1)Kr^{-2n}\frac{1}{2}r^2 d\eta.
\end{equation}
Therefore we have
\begin{equation}
 \begin{split}
  (\omega+dd^c \rho)^n &= n(1+2(n-1)^2 Kr^{-2n})(1-2(n-1)Kr^{-2n})^{n-1} rdr\wedge\eta\wedge(\frac{1}{2}r^2 d\eta)^{n-1} \\
               &= (1+2(n-1)^2 Kr^{-2n})(1-2(n-1)Kr^{-2n})^{n-1} \omega^n \\
            & \leq\omega^n,
 \end{split}
\end{equation}
for $K>0$ and $2K(n-1)r^{-2n}\leq 1$.  Then for suitably large $K>0$, with $r_0$ possibly increased, one has $\phi\leq\rho$ on $T_{r_0}$.
An application of the maximal principle using that $\phi\rightarrow 0$ as $y\rightarrow\infty$ gives the upper bound in (\ref{eq:conv}).

For the lower bound set $\rho =Kr^{-2n+2}(\log r)^\delta$.  Then a similar computation gives
\begin{equation}
\begin{split}
 (\omega +dd^c\rho)^n & =\left(1-\delta(n-1)Kr^{-2n}(\log r)^{\delta-1} +o(r^{-2n}(\log r)^{\delta-1})\right)\omega^n\\
                      & \geq\omega^{n},
\end{split}
\end{equation}
on $Y_{r_0}$ for $K<0$ and $r_0>0$ sufficiently large.  And another application of the maximum principle give the lower bound in (\ref{eq:conv}).
\end{proof}

The following proposition is a slight variation of proposition 5.1 in~\cite{TY2}.  We give a somewhat simpler proof
for this context.
\begin{prop}\label{prop:conv2}
Let $\phi$ be as above.  Then for $\frac{1}{2}>\delta>0$, there are constants $C_{\delta,k}$ depending only on
$k$ and $\delta$ so that
\begin{equation}
\|\nabla^k \phi\|_{g_0}(y) \leq C_{\delta,k}r(y)^{-2n+2-k+\delta},\quad\text{for },y\in Y_{r_0}.
\end{equation}
\end{prop}
\begin{proof}
Recall that on $Y_{r_0}$ $g_0 =dr^2 +r^2 g$, the cone metric, and the Euler vector field $r\partial_r$
generates an action of $\R_{>0}$ by homothetic isometries on $g_0$.  For $a>1$ denote this action by
$\psi_a :Y_{r_0}\rightarrow Y_{r_0}$.  Then
\begin{equation}\label{eq:homo}
\psi_a^* g_0 =a^2 g_0.
\end{equation}
Then it easily follows that for all $k\geq 0$
\begin{equation}\label{eq:curv-cone}
\|\nabla^k R(g_0)\|_{g_0} =O(r^{-k-2}).
\end{equation}

Let $b>r_0$, and $S_b \subset Y_{r_0}$ be the link $S_b =\{y\in Y_{r_0}: r(y)=b\}$.  Cover $S_b$ with
coordinate balls $B_\rho (x,g_0), x\in S_b$ of radius $\rho$ so that $B_{\frac{\rho}{2}}(x,g_0)$ cover $S_b$.  Then for $a>1$
$\psi_a (B_\rho (x,g_0))=B_{a\rho}(\psi_a(x),g_0)$.  Define $\phi_a :=a^{-2}\psi_a^*\phi$.
Since
\begin{equation}
\psi_a^* \|\nabla^k \phi\|_{g_0} =a^{-k}\|\nabla^k \psi_a^*\phi\|_{g_0} =a^{-k+2}\|\nabla^k \phi_a\|_{g_0},
\end{equation}
it is sufficient to show there are constants $C_{\delta,k}$ so that
\begin{equation}\label{eq:suffic}
\|\nabla^k \phi_a\|_{g_0}(x) \leq C_{\delta,k}a^{-2n+\delta},\quad\text{for }x\in S_b.
\end{equation}
Note that (\ref{eq:suffic}) holds for $k=0$ by Lemma~\ref{lem:conv1}.

Recall that, on $Y_{r_0}$, $\phi$ is a solution to
\begin{equation}\label{eq:monge-loc1}
(\omega_0 +i\partial\ol{\partial}\phi)^n =\omega_0^n.
\end{equation}
Apply $\psi_a^*$ to (\ref{eq:monge-loc1}) and rescale to get
\begin{equation}\label{eq:monge-loc2}
(\omega_0 +i\partial\ol{\partial}\phi_a)^n =\omega_0^n.
\end{equation}

Let $\omega=\omega_0 +i\partial\ol{\partial}\phi$, and define an operator $P$ on $B_\rho (x,g_0), x\in S_b$
by
\begin{equation}\label{eq:scha1}
(Pu)\omega_0^n := i\partial\ol{\partial}u\wedge(\omega^{n-1} +\omega^{n-2}\omega_0 +\cdots +\omega_0^{n-1}).
\end{equation}
Then since the proof of Proposition~\ref{prop:monge} gives a bound on $\phi$ in $C^{2,\alpha},0< \alpha<1$
and $c^{-1}\omega_0 \leq\omega\leq c\omega_0$ for some $c>0$, the Schauder interior estimates
(cf.~\cite{GiTr} Theorem 6.2) apply to (\ref{eq:scha1}).  Thus if $Pu =f$ with $u\in C^2 (B_\rho)$ and
$f\in C^{0,\alpha}(B_\rho)$, then $u|_{B_{\frac{\rho}{2}}}\in C^{2,\alpha}(B_{\frac{\rho}{2}})$ and
\begin{equation}\label{eq:schauder}
\|u|_{B_{\frac{\rho}{2}}}\|_{C^{2,\alpha}}\leq C\left(\|f\|_{C^{0,\alpha}} + \|u\|_{C^0}\right).
\end{equation}

Then from (\ref{eq:monge-loc2}) we have $P(\phi_a)=0$.  So
\begin{equation}
\|\phi_a |_{B_{\frac{\rho}{2}}}\|_{C^{2,\alpha}} \leq C\|\phi_a |_{B_\rho}\|_{C^{0}}\leq C' a^{-2n+\delta}.
\end{equation}

Now apply the covariant derivative $\nabla_e$ to (\ref{eq:monge-loc2}) to get
\begin{equation}\label{eq:scha2}
g^{\beta\ol{\gamma}}\nabla_\beta \nabla_{\ol{\gamma}}\nabla_e \phi_a =R^{\ol{\epsilon}}_{\ol{\gamma}e\beta}\nabla_{\ol{\epsilon}}\phi_a,
\end{equation}
where $g$ is the metric with K\"{a}hler form $\omega$, while the covariant derivative and curvature is with respect
to $g_0$.  The $C^{0,\alpha}$ norm of the right-hand side of
(\ref{eq:scha2}) is bounded by $Ca^{-2n+\delta}$ for some $C$.  So (\ref{eq:schauder}) implies that
\begin{equation}
\|\phi_a |_{B_{\frac{\rho}{2}}}\|_{C^{3,\alpha}}\leq Ca^{-2n+\delta}.
\end{equation}

We proceed inductively.  Suppose we have the bound
\begin{equation}
\|\phi_a |_{B_{\rho}}\|_{C^{k,\alpha}}\leq Ca^{-2n+\delta}.
\end{equation}
Apply the general $k-1=i+j$ order covariant derivative
$\nabla_{\alpha_1}\cdots\nabla_{\alpha_i}\nabla_{\ol{\epsilon}_1}\cdots\nabla_{\ol{\epsilon_j}}$ to (\ref{eq:monge-loc2})
and rearrange terms using curvature identities to get
\begin{equation}\label{eq:scha-ind}
g^{\beta\ol{\gamma}}\nabla_\beta\nabla_{\ol{\gamma}}\nabla_{\alpha_1}\cdots\nabla_{\alpha_i}\nabla_{\ol{\epsilon}_1}\cdots\nabla_{\ol{\epsilon_j}}
\phi_a =F,
\end{equation}
where $F$ is an expression containing the curvature tensor, its covariant derivatives, and the covariant derivatives of $\phi_a$ up to order $k-1$.  Thus $F$ is bounded in $C^{1,\alpha}$ by $Ca^{-2n+\delta}$ on $B_{\rho}$ by the previous step.  Then apply (\ref{eq:schauder})
to the equation (\ref{eq:scha-ind}) to get a bound
\begin{equation}
\|\phi_a |_{B_{\frac{\rho}{2}}}\|_{C^{k+1,\alpha}}\leq Ca^{-2n+\delta}.
\end{equation}
\end{proof}

Theorem~\ref{thm:main} now follows from Proposition~\ref{prop:monge}, Lemma~\ref{lem:conv1}, and
Proposition~\ref{prop:conv2}.

We now collect some of the asymptotic properties of the metric in Theorem~\ref{thm:main} which follow from the
preceding results and equation (\ref{eq:curv-cone}).
\begin{prop}\label{prop:decay}
Let $g$ be the Ricci-flat K\"{a}hler metric on $Y$ of Proposition~\ref{prop:monge}.  Then
curvature of $g$ satisfies
\begin{equation}
\|\nabla^k R(g)\|_{g} =O(r^{-2-k}), \quad\text{ for } k\geq 0.
\end{equation}

Furthermore, if $\|R(g)\|_g =O(r^{-\alpha})$, for $\alpha >2$, then $(Y,g)$ is asymptotically locally Euclidean
of order $2n$.
\end{prop}

The second statement of the Proposition follows from a result of~\cite{BKN}.  We recall the definition
of asymptotically locally Euclidean (ALE).
By ALE of order $m$ we mean the following.  There exists a compact subset $K\subset Y$, a finite group
$\Gamma\subset O(2n)$ acting freely on $\R^{2n}\setminus\{0\}$, and a ball $B_R(0)\subset\R^{2n}$ of radius $R>0$.
So that there is a diffeomorphism $\chi:\R^{2n}/\Gamma \rightarrow Y\setminus K$ and
\begin{equation}
 \nabla^k \chi^*g -\nabla^k h=O(r^{-m-k}),
\end{equation}
where $h$ is the flat metric and $\nabla$ its covariant derivative.

Furthermore, since $Y$ is K\"{a}hler it is not difficult to show that one may take $\R^{2n}=\C^n$ with the standard
complex structure $J_0$ and $\Gamma\subset U(n)$.  And if $J$ is the complex structure on $Y$ we have
\begin{equation}
 \nabla^k \chi^*J -\nabla^k J_0=O(r^{-m-k}),
\end{equation}
and Ricci-flatness implies that $\Gamma\subset SU(n)$.  The results of~\cite{BKN} imply that if
$\|R(g)\|_g =O(r^{-\alpha})$, for $\alpha >2$, then $(Y,g)$ is ALE of order $2n$.

\section{Toric case}

We now restrict to the toric case.  We will consider crepant resolutions $\pi:Y\rightarrow X$ where both $X$ and $Y$ are
toric varieties.  In this case $X=C(S)$ is a toric K\"{a}hler cone over a toric Sasaki manifold $S$.
We will prove the toric version of Theorem~\ref{thm:main}, Corollary~\ref{cor:main}, which makes use
of the general existence result of A. Futaki, H. Ono, and G. Wang~\cite{FOW} of Ricci-flat K\"{a}hler cone metrics on
$X=C(S)$ provided $S$ satisfies the condition in Proposition~\ref{prop:CY-cond-toric},
which is a translation into toric geometry of the condition in Proposition~\ref{prop:CY-cond}.
Then it is elementary using toric geometry to construct examples
of crepant resolutions $\pi:Y\rightarrow X$ of a Ricci-flat K\"{a}hler cone $X$.
We will start with the differential geometric picture of toric geometry.  See~\cite{Gui2} for a good reference.  Then
we will use concepts from the algebraic geometric picture of toric varieties to construct crepant resolutions.
A good reference for this is~\cite{Od}.

\subsection{Toric Sasaki-Einstein manifolds}

In this section we recall the basics of toric Sasaki manifolds.  Much of what follows can be found in
~\cite{MSY} or~\cite{FOW}.

\begin{defn}
 A Sasaki manifold $(S,g)$ of dimension $2n-1$ is \emph{toric} if there is an effective action of an
 $n$-dimensional torus $T=T^{n}$ preserving the Sasaki structure such that the Reeb vector field $\xi$ is an
 element of the Lie algebra $\mathfrak{t}$ of $T$.

 Equivalently, a toric Sasaki manifold is a Sasaki manifold $S$ whose K\"{a}hler cone
 $C(S)$ is a toric K\"{a}hler manifold.
\end{defn}

We have an effective holomorphic action of $T_\C \cong (\C^*)^{n}$ on $C(S)$ whose restriction to
$T\subset T_\C$ preserves the K\"{a}hler form $\omega =d(\frac{1}{2}r^2 \eta)$.  So there is a moment map
\begin{equation}\label{eq:moment-map}
\begin{gathered}
 \mu: C(S) \longrightarrow \mathfrak{t}^* \\
 \langle \mu(x),X\rangle = \frac{1}{2}r^2\eta(X_S (x)),
\end{gathered}
\end{equation}
where $X_S$ denotes the vector field on $C(S)$ induced by $X\in\mathfrak{t}$.  We have the
moment cone defined by
\begin{equation}
 \mathcal{C}(\mu) :=\mu(C(S)) \cup \{o\},
\end{equation}
which from~\cite{Ler} is a strictly convex rational polyhedral cone.  Recall that this means that there are vectors
$u_i,i=1,\ldots,d$, in the integral lattice $\Z_T =\ker\{\exp(2\pi\cdot):\mathfrak{t}\rightarrow T\}\subset\mathfrak{t}$ such that
\begin{equation}\label{eq:moment-cone}
 \mathcal{C}(\mu) =\bigcap_{j=1}^{d} \{y\in\mathfrak{t}^* : \langle u_j,y\rangle\geq 0\}.
\end{equation}
The condition that $\mathcal{C}(\mu)$ is strictly convex means that it is not contained in any linear subspace of $\mathfrak{t}^*$,
and it is cone over a finite polytope.
We assume that the set of vectors $\{u_j\}$ is minimal in that removing one changes the set defined by
(\ref{eq:moment-cone}).  And we furthermore assume that the vectors $u_j$ are primitive, meaning that
$u_j$ cannot be written as $p\tilde{u}_j$ for $p\in\Z,p>1,$ and $\tilde{u}_j\in\Z_T$.

Let $\inter\mathcal{C}(\mu)$ denote the interior of $\mathcal{C}(\mu)$.  Then the action of $T$ on
$\mu^{-1}(\inter\mathcal{C}(\mu))$ is free and it is a Lagrangian torus fibration over $\inter\mathcal{C}(\mu)$.
There is a condition on the $\{u_j\}$ for $S$ to be a smooth manifold.  Each face
$\mathcal{F}\subset\mathcal{C}(\mu)$ is the intersection of a number of facets
$\{y\in\mathfrak{t}^* :l_j(y)=\langle u_j ,y \rangle =0\}$.  Let $u_{j_1},\ldots,u_{j_a}$ be the corresponding
collection of normal vectors in $\{u_j\}$, where $a$ is the codimension of $\mathcal{F}$.  Then $S$ is smooth, and
the cone $\mathcal{C}(\mu)$ is said to be \emph{non-singular} if and only if
\begin{equation}\label{eq:nonsing}
 \left\{ \sum_{k=1}^{a} \nu_{k} u_{j_k} :\nu_k \in\R\right\}\cap\Z_T =\left\{\sum_{k=1}^{a} \nu_{k} u_{j_k} :\nu_k \in\Z\right\}
\end{equation}
for all faces $\mathcal{F}$.

Note that $\mu(S)=\{y\in\mathcal{C}(\mu) : y(\xi)=\frac{1}{2}\}$.  The hyperplane
$\{y\in\mathfrak{t}^* : y(\xi)=\frac{1}{2}\}$ is called the \emph{characteristic hyperplane} of the Sasaki structure.
Consider the dual cone to $\mathcal{C}(\mu)$
\begin{equation}\label{eq:dual-cone}
 \mathcal{C}(\mu)^* =\{\tilde{x}\in\mathfrak{t} :\langle\tilde{x}, y\rangle\geq 0\text{ for all }y\in\mathcal{C}(\mu)\},
\end{equation}
which is also a strictly convex rational polyhedral cone by Farkas' theorem.  Then $\xi$ is in the interior of
$\mathcal{C}(\mu)^*$.  Let $\frac{\partial}{\partial\phi_i},i=1,\ldots, n$ be a basis of $\mathfrak{t}$ in $\Z_T$.
Then we have the identification $\mathfrak{t}^*\cong\mathfrak{t}\cong\R^{n}$ and we write
\[ u_j =(u_j^1, \ldots,u_j^{n}),\quad \xi=(\xi^1,\ldots,\xi^{n}). \]
If we set
\begin{equation}\label{eq:symp-coord}
y_i =\langle\mu(x),\frac{\partial}{\partial\phi_i}\rangle\quad ,i=1,\ldots, n,
\end{equation}
then we have symplectic coordinates $(y,\phi)$ on $\mu^{-1}(\inter\mathcal{C}(\mu))\cong \inter\mathcal{C}(\mu)\times T^{n}$.
In these coordinates the symplectic form is
\begin{equation}
 \omega =\sum_{i=1}^{n} dy_i \wedge d\phi_i.
\end{equation}
The K\"{a}hler metric can be seen as in~\cite{Abr} to be of the form
\begin{equation}
 g=\sum_{ij} G_{ij}dy_i dy_j + G^{ij}d\phi_i d\phi_j,
\end{equation}
where $G^{ij}$ is the inverse matrix to $G_{ij}(y)$, and the complex structure is
\begin{equation}\label{eq:comp-st}
 \mathcal{I} =\left\lgroup
\begin{matrix}
 0 & -G^{ij} \\
 G_{ij} & 0 \\
\end{matrix}\right\rgroup
\end{equation}
in the coordinates $(y,\phi)$.  The integrability condition of $\mathcal{I}$ is equivalent to $G_{ij,k} =G_{ik,j}$.  Thus
\begin{equation}
 G_{ij}=G_{,ij} :=\frac{\partial^2 G}{\partial y_i \partial y_j},
\end{equation}
for some strictly convex function $G(y)$ on $\inter\mathcal{C}(\mu)$.  We call $G$ the symplectic potential of the
K\"{a}hler metric.

One can construct a canonical K\"{a}hler structure on the cone $X=C(S)\cup\{o\}$, with a fixed holomorphic structure,
via a simple K\"{a}hler reduction of $\C^d$ (cf.~\cite{Gui1,Gui2} and~\cite{BurGuiLer}).  This procedure will be recounted
in Section~\ref{subsect:crep-res}.

The symplectic potential of the canonical K\"{a}hler metric is
\begin{equation}
 G^{can} =\frac{1}{2}\sum_{i=1}^{d}l_i (y)\log l_i (y).
\end{equation}
Let
\[ G_\xi =\frac{1}{2}l_{\xi}(y)\log l_{\xi} -\frac{1}{2}l_{\infty}(y)\log l_{\infty}(y),\]
where
\[ l_{\xi}(y)=\langle\xi, y\rangle, \text{ and } l_{\infty}(y)=\sum_{i=1}^d \langle u_i ,y\rangle.\]
Then
\begin{equation}\label{eq:can-pot}
 G_\xi ^{can} =G^{can} + G_\xi,
\end{equation}
defines a symplectic potential of a K\"{a}hler metric on $C(S)$ with induced Reeb vector field $\xi$.
To see this write
\begin{equation}
 \xi =\sum_{i=1}^{n} \xi^i \frac{\partial}{\partial\phi_i},
\end{equation}
and note that the Euler vector field is
\begin{equation}
 r\frac{\partial}{\partial r} =2\sum_{i=1}^{n}y_i\frac{\partial}{\partial y_i}.
\end{equation}
Thus from (\ref{eq:comp-st}) we must have
\begin{equation}\label{eq:sym-reeb}
 \xi^i =\sum_{j=1}^{n} 2G_{ij}y_j.
\end{equation}
Computing from (\ref{eq:can-pot}),
\begin{equation}\label{eq:can-pot-der}
 \left(G_\xi ^{can} \right)_{ij} =\frac{1}{2}\sum_{k=1}^d \frac{u_k^i u_k^j}{l_k (y)} +\frac{1}{2}\frac{\xi^i \xi^j}{l_\xi (y)}
 -\frac{1}{2}\frac{\sum_{k=1}^d u_k^i \sum_{k=1}^d u_k^j}{l_\infty (y)}.
\end{equation}
And plugging (\ref{eq:can-pot-der}) into (\ref{eq:sym-reeb}) shows we have the desired Reeb vector field.

The general symplectic potential is of the form
\begin{equation}
 G =G^{can} + G_\xi +g,
\end{equation}
where $g$ is a smooth homogeneous degree one function on $\mathcal{C}$ such that $G$ is strictly convex.
The following follows easily from this discussion.
\begin{prop}
 Let $S$ be a compact toric Sasaki manifold and $C(S)$ its K\"{a}hler cone.  For any
 $\xi\in\inter\mathcal{C}(\mu)^*$ there exists a toric K\"{a}hler cone metric, and associated Sasaki structure on $S$,
 with Reeb vector field $\xi$.  And any other such structure is a transverse K\"{a}hler deformation, i.e.
 $\tilde{\eta} =\eta +2d^c \phi$, for a $T$-invariant function $\phi$.
\end{prop}

We consider now the holomorphic picture of $C(S)$.  Note that the complex structure on $X=C(S)$ is determined up to
biholomorphism by the associated moment polyhedral cone $\mathcal{C}(\mu)$ (cf.~\cite{Abr} Proposition A.1).
And the construction of $X=C(S)$ as in~\cite{Gui1,Gui2} shows that $X=C(S)$ is a toric variety
with open dense orbit $(\C^*)^{n}\cong\mu^{-1}(\inter\mathcal{C})\subset C(S)$.

Recall that a toric variety is characterized by a \emph{fan} (cf.~\cite{Od}).  We give some definitions.
\begin{defn}
A subset $\sigma$ of $\mathfrak{t}\cong\R^n$ is a \emph{strongly convex rational polyhedral cone}, if
there exists a finite number of elements $u_1,u_2,\ldots,u_s$ in $\Z_T \cong\Z^n$
such that
\[\sigma =\{a_1 u_1 +\cdots+a_s u_s :a_i \in\R_{\geq 0}\text{ for }i=1,\ldots,s\},\]
and $\sigma\cap(-\sigma)=\{o\}$.
\end{defn}
\begin{defn}
 A \emph{fan} in $\Z_T \cong\Z^n$ is a nonempty collection $\Delta$  of strongly convex rational polyhedral cones
 in $\mathfrak{t}\cong\R^n$ satisfying the following:
\begin{thmlist}
\item  Every face of any $\sigma\in\Delta$ is contained in $\Delta$.

\item  For any $\sigma,\sigma'\in\Delta$, the intersection $\sigma\cap\sigma'$ is a face of both $\sigma$ and $\sigma'$.
\end{thmlist}
\end{defn}

Then to every fan $\Delta$ in $\Z_T \cong\Z^n$ is uniquely associated a normal complex algebraic variety $X_{\Delta}$
with an algebraic action of $T_{\C}\cong(\C^*)^n$.  Furthermore, there is an open dense orbit isomorphic to
$T_{\C}\cong(\C^*)^n$.  Conversely, if a torus $(\C^*)^n$ acts algebraically on a normal algebraic variety $X$, with
locally finite type over $\C$, with an open dense orbit isomorphic to $(\C^*)^n$, then there is a fan $\Delta$ in
$\Z^n$ with $X$ equivariantly isomorphic to $X_{\Delta}$.  See~\cite{Od} for more details.

There is a fan in $\Z_T \subset\mathfrak{t}$ associated to every strictly convex rational polyhedral set
$\mathcal{C}\subset\mathfrak{t}^*$.  Suppose
\begin{equation}
\mathcal{C}=\bigcap_{j=1}^{d} \{y\in\mathfrak{t}^* : \langle u_j,y\rangle\geq\lambda_j\},
\end{equation}
where $u_j \in\Z_T$ and $\lambda_j \in\R$ for $j=0,\ldots,d$.  Each face $\mathcal{F}\subset\mathcal{C}$ is the
intersection of facets $\{y\in\mathfrak{t}^* : l_{j_k}(y)=\langle u_{j_k},y\rangle-\lambda_{j_k} =0\}\cap\mathcal{C}$
for $k=1,\ldots, a$, where $\{j_1,\ldots,j_a\}\subseteq\{1,\ldots,d\}$, and the codimension of $\mathcal{F}$ is $a$.
Then to the face $\mathcal{F}$ we associate a cone $\sigma_{\mathcal{F}}$ in $\mathfrak{t}\cong\R^n$
\begin{equation}\label{eq:dual-cone-face}
\sigma_{\mathcal{F}} =\{c_1 u_{j_1} +\cdots+c_a u_{j_a} :c_k \in\R_{\geq 0}\text{ for }k=1,\ldots,a\}.
\end{equation}
It is easy to see that the set of all $\sigma_{\mathcal{F}}$ for faces $\mathcal{F}\subseteq\mathcal{C}$ define
a fan $\Delta$ in $\Z_T$.

Consider the convex polyhedral cone $\mathcal{C}(\mu)$.  From (\ref{eq:moment-cone}) the fan in $\Z_T$ associated to
$\mathcal{C}(\mu)$ consists of the dual cone (\ref{eq:dual-cone}) and all of its faces, where the dual cone is
\begin{equation}\label{eq:dual-cone-fan}
\mathcal{C}(\mu)^* =\{c_{1} u_1 +\cdots +c_d u_d : c_k\in\R_{\geq 0}\text{ for }k=1,\ldots,d\}.
\end{equation}
It follows that $C(S)$ is an affine variety as its fan has a single n-dimensional cone.

We introduce logarithmic
coordinates $(z_1,\ldots,z_{n}) =(x_1 +i\phi_1,\ldots, x_{n}+i\phi_{n})$ on
$\C^{n}/{2\pi i\Z^{n}}\cong (\C^*)^{n}\cong\mu^{-1}(\inter\mathcal{C})\subset C(S)$, i.e.
$x_j +i\phi_j =\log w_j$ if $w_j,j=1,\ldots,n$, are the usual coordinates on $(\C^*)^{n}$.
Since on $\mu^{-1}(\inter\mathcal{C})$ the K\"{a}hler form $\omega$ is $T^n$ invariant and the $T^n$-action is
Hamiltonian, we have
\begin{equation}
 \omega =\mathbf{i}\partial\ol{\partial}F,
\end{equation}
where $F$ is a strictly convex function of $(x_1,\ldots,x_{n})$ (cf.~\cite{Gui1} Theorem 4.3).  One can check that
\begin{equation}
 F_{ij}(x) =G^{ij}(y),
\end{equation}
where $\mu=y=\frac{\partial F}{\partial x}$ is the moment map.  Strictly speaking, $\mu=\frac{\partial F}{\partial x} +c$ for a constant
$c\in\R^n$.  But we add a linear factor to $F$, so that $\mu=\frac{\partial F}{\partial x}$.
Furthermore, one can show $x=\frac{\partial G}{\partial y}$, and
the K\"{a}hler and symplectic potentials are related by the Legendre transform
\begin{equation}\label{eq:legen}
 F(x)= \sum_{i=1}^{n} x_i \cdot y_i -G(y).
\end{equation}
It follows from equation (\ref{eq:symp-coord}) defining symplectic coordinates that
\begin{equation}\label{eq:pot-id}
 F(x)= l_\xi (y) =\frac{r^2}{2}.
\end{equation}
The potential $F(x)$ is, of course, only defined up to an affine function on $(x_1,\ldots,x_{n})$, but by considering the
limit as $y\mapsto 0$ in (\ref{eq:legen}) one shows that the first equality in (\ref{eq:pot-id}) holds.

We now consider the conditions in Proposition~\ref{prop:CY-cond} more closely in the toric case.
So suppose the Sasaki structure satisfies Proposition~\ref{prop:CY-cond}, thus we may assume
$c_1^B =2n[\omega^T]$.  Then equation (\ref{eq:ricci-cone}) implies that
\begin{equation}
 \rho =-i\partial\ol{\partial}\log\det(F_{ij})=i\partial\ol{\partial}h,
\end{equation}
with $\xi h=0=r\frac{\partial}{\partial r}h$, and we may assume $h$ is $T^{n}$-invariant.
Since a $T^{n}$-invariant pluriharmonic function is an affine function, we have
constants $\gamma_1,\ldots,\gamma_{n}\in\R$ so that
\begin{equation}\label{eq:hol-CY}
 \log\det(F_{ij})=-2\sum_{i=1}^{n}\gamma_i x_i -h.
\end{equation}
In symplectic coordinates we have
\begin{equation}\label{eq:symp-CY}
 \det(G_{ij})=\exp(2\sum_{i=1}^{n}\gamma_i G_i +h).
\end{equation}
Then from (\ref{eq:can-pot}) one computes the right hand side to get
\begin{equation}
 \det(G_{ij})=\prod_{k=1}^d \left(\frac{l_k (y)}{l_\infty (y)}\right)^{(\gamma,u_k)} (l_\xi (y))^{-n} \exp(h),
\end{equation}
And from (\ref{eq:can-pot-der}) we compute the left hand side of (\ref{eq:symp-CY})
\begin{equation}
 \det(G_{ij})=\prod_{k=1}^d(l_k (y))^{-1} f(y),
\end{equation}
where $f$ is a smooth function on $\mathcal{C}(\mu)$.  Thus $(\gamma, u_k)=-1$, for $k=1,\ldots,d$.
Since $\mathcal{C}(\mu)^*$ is strictly convex, $\gamma$ is a uniquely determined element of $\mathfrak{t}^*$.

Applying $\sum_{j=1}^{n}y_j \frac{\partial}{\partial y_j}$ to (\ref{eq:symp-CY}) and noting that
$\det(G_{ij})$ is homogeneous of degree $-n$ we get
\begin{equation}\label{eq:reeb-const}
 (\gamma,\xi)=-n.
\end{equation}

As in Proposition~\ref{prop:CY-cond} $e^h \det(F_{ij})$ defines a flat metric $\|\cdot\|$ on
$\mathbf{K}_{C(S)}$.  Consider the $(n,0)$-form
\[ \Omega = e^{i\theta}e^{\frac{h}{2}}\det(F_{ij})^{\frac{1}{2}} dz_1 \wedge\cdots\wedge dz_{n}.\]
From equation (\ref{eq:hol-CY}) we have
\[ \Omega =e^{i\theta}\exp(-\sum_{j=1}^{n}\gamma_j x_j)dz_1 \wedge\cdots\wedge dz_{n}.\]
If we set $\theta=-\sum_{j=1}^{n}\gamma_j \phi_j$, then
\begin{equation}\label{eq:hol-form-toric}
 \Omega=e^{-\sum_{j=1}^{n}\gamma_j z_j}dz_1 \wedge\cdots\wedge dz_{n}
\end{equation}
is clearly holomorphic on $U=\mu^{-1}(\inter\mathcal{C})$.  When $\gamma$ is not integral, then we take $\ell\in\Z_+$
such that $\ell\gamma$ is a primitive element of $\Z_T^* \cong\Z^{n}$.  Then
$\Omega^{\otimes\ell}$ is a holomorphic section of $\mathbf{K}^{\ell}_{C(S)}|_U$ which extends to
a holomorphic section of $\mathbf{K}^{\ell}_{C(S)}$ as $\|\Omega\|=1$.

It follows from (\ref{eq:hol-form-toric}) that
\begin{equation}
 \mathcal{L}_\xi \Omega =-i(\gamma,\xi)\Omega =i n\Omega.
\end{equation}
And note that we have equation (\ref{eq:hol-form}) from (\ref{eq:hol-CY}) and (\ref{eq:hol-form-toric}).
We collect these results in the following proposition.

\begin{prop}\label{prop:CY-cond-toric}
Let $S$ be a compact toric Sasaki manifold of dimension $2n-1$.  Then the conditions of
Proposition~\ref{prop:CY-cond} are equivalent to the existence of $\gamma\in\mathfrak{t}^*$ such that
\begin{thmlist}
 \item $(\gamma, u_k)=-1$, for $k=1,\ldots,d$,

 \item $(\gamma,\xi)=-n$, and

 \item there exists $\ell\in\Z_+$ such that $\ell\gamma\in\Z_T^* \cong\Z^{n}$.
\end{thmlist}
Then (\ref{eq:hol-form-toric}) defines a nowhere vanishing section of $\mathbf{K}^{\ell}_{C(S)}$.  And
$C(S)$ is $\ell$-Gorenstein if and only if a $\gamma$ satisfying the above exists.
\end{prop}

We will need the beautiful result of A. Futaki, H. Ono, and G. Wang on the existence of
Sasaki-Einstein metrics on toric Sasaki manifolds.
\begin{thm}[\cite{FOW,CFO}]\label{thm:FOW}
Suppose $S$ is a toric Sasaki manifold satisfying Proposition~\ref{prop:CY-cond-toric}.
Then we can deform the Sasaki structure by varying the Reeb vector field and then performing
a transverse K\"{a}hler deformation to a Sasaki-Einstein metric.  The Reeb vector field and transverse
K\"{a}hler deformation are unique up to isomorphism.
\end{thm}

In~\cite{FOW} a more general result is proved.  It is proved that a compact toric Sasaki
manifold satisfying Proposition~\ref{prop:CY-cond-toric} has a transverse K\"{a}hler deformation
to a Sasaki structure satisfying the transverse K\"{a}hler Ricci soliton equation:
\[ \rho^T - 2n\omega^T =\mathcal{L}_X \omega^T \]
for some Hamiltonian holomorphic vector field $X$.  The analogous result for toric Fano manifolds was
proved in~\cite{WaZh}.  A transverse K\"{a}hler Ricci soliton becomes a transverse K\"{a}hler-Einstein
metric, i.e. $X=0$, if the Futaki invariant $f_1$ of the transverse K\"{a}hler structure vanishes.
The invariant $f_1$ depends only on the Reeb vector field $\xi$.  The next step is to use a
volume minimization argument due to Martelli-Sparks-Yau~\cite{MSY} to show there is a unique $\xi$
satisfying (\ref{eq:reeb-const}) for which $f_1$ vanishes.

\begin{xpl}\label{xpl:two-points}
Let $M=\cps^2_{(2)}$ be the two-points blow up.  And Let $S \subset\mathbf{K}_M$ be the
$U(1)$-subbundle of the canonical bundle.  Then the standard Sasaki structure on $S$
satisfies (i) of Proposition~\ref{prop:CY-cond}, and it is not difficult to show that $S$
is simply connected and is toric.
But the automorphism group of $M$ is not reductive, thus $M$ does not admit a K\"{a}hler-Einstein metric
due to Y. Matsushima~\cite{Mat}.  Thus there is no Sasaki-Einstein structure with the usual Reeb vector field.
But by Theorem~\ref{thm:FOW} there is a Sasaki-Einstein structure with a different Reeb vector field.

The vectors defining the facets of $\mathcal{C}(\mu)$ are
\[u_1 =(0,0,1), u_2 =(0,1,1) ,u_3 =(1,2,1) ,u_4 =(2,1,1) ,u_5 =(1,0,1).\]
The Reeb vector field of the toric Sasaki-Einstein metric on $S$ was calculated in~\cite{MSY} to be
\[\xi =\left(\frac{9}{16}(-1+\sqrt{33}), \frac{9}{16}(-1+\sqrt{33}), 3\right).\]
One sees that the Sasaki structure is irregular with the closure of the generic orbit being a two torus.
\end{xpl}

\subsection{Toric crepant resolutions}\label{subsect:crep-res}

Let $X=C(S)$ be a toric K\"{a}hler cone.  Then as an algebraic variety $X=X_{\Delta}$ where $\Delta$ is the fan in
$\Z_T \cong\Z^n$ defined by the dual cone $\mathcal{C}(\mu)^*$, spanned by $u_1,\ldots,u_d\in\Z_T$, and its faces as
in (\ref{eq:dual-cone-fan}).  We assume that $X$ is Gorenstein.  Thus there is a $\gamma\in\Z_T^*$ so that
$\gamma(u_i)=-1$ for $i=1,\ldots, d$.  Let $H_\gamma =\{x\in\mathfrak{t}: \langle\gamma,x\rangle =-1\}$ be the hyperplane
defined by $\gamma$.  Then
\begin{equation}
P_{\Delta}:=\{x\in \mathcal{C}(\mu)^* : \langle\gamma,x\rangle =-1\}\subset H_\gamma \cong\R^{n-1}
\end{equation}
is an $(n-1)$-dimensional lattice polytope.  The lattice being $H_\gamma \cap\Z_T \cong\Z^{n-1}$.

A toric crepant resolution
\begin{equation}\label{eq:toric-resol}
\pi: X_{\tilde{\Delta}}\rightarrow X_{\Delta}
\end{equation}
is given by a nonsingular subdivision
$\tilde{\Delta}$ of $\Delta$ with every 1-dimensional cone $\tau_i \in\tilde{\Delta}(1), i=1,\ldots,N$, generated
by a primitive vector $u_i :=\tau_i \cap H_\gamma$.  This is equivalent to a basic, lattice triangulation of
$P_{\Delta}$.  \emph{Lattice} means that the vertices of every simplex are lattice points, and \emph{basic}
means that the vertices of every top dimensional simplex generates a basis of $\Z^{n-1}$.
Note that a \emph{maximal} triangulation of $P_{\Delta}$, meaning that the vertices of every simplex are its only
lattice points, always exists.  Every basic lattice triangulation is maximal, but the converse only holds in
dimension 2.  In dimensions $\geq 3$ there are polytopes which do not admit basic lattice triangulations.

The condition that $u_i :=\tau_i \cap H_\gamma$ is primitive for each $i=1,\ldots,N$ is precisely the condition
that the section of Proposition (\ref{prop:CY-cond-toric}), $\Omega\in\Gamma(\mathbf{K}_{C(S)})$, characterized
by $\gamma\in\Z_T$ lifts to a non-vanishing section of $\mathbf{K}_{X_{\tilde{\Delta}}}$.
See~\cite{Od}, Proposition 2.1.

Note that a toric crepant resolution (\ref{eq:toric-resol}) of $X_{\Delta}$ is not unique, if one exists.
But if $E=\pi^{-1}(o)$ is the exceptional set, then the number of prime divisors in $E$ is invariant.
There is a prime divisor $D_i, i=d+1,\ldots,N$, for each lattice point in $\inter P_{\Delta}$.

\begin{prop}\label{prop:3dim-resol}
Suppose $X=C(S)$ is a 3-dimensional, $n=3$, Gorenstein toric K\"{a}hler cone.  Then $P_{\Delta}$ admits a
basic lattice triangulation.  Thus $X$ admits a toric crepant resolution.
\end{prop}
\begin{proof}
There is a maximal lattice triangulation of $P_{\Delta}$. Since it is 2-dimensional, any maximal triangulation
is basic.
\end{proof}

One can further show that any 3-dimensional Gorenstein toric K\"{a}hler cone admits a crepant resolution satisfying
the requirements of Corollary~\ref{cor:main}.  See~\cite{vC3}.

In the previous section we associated a fan in $\Z_T$ to every rational convex polyhedral set
$\mathcal{C}\subset\mathfrak{t}^*$.  The following definition will be used to associate a rational convex
polyhedral set to a fan.
\begin{defn}
A real valued function $h:|\Delta|\rightarrow\R$ on the support $|\Delta|:=\cup_{\sigma\in\Delta}\sigma$ is a
\emph{support function} if it is linear on each $\sigma\in\Delta$.  That is, there exist an
$l_{\sigma}\in(\R^n)^*$ for each $\sigma\in\Delta$ so that $h(x)=\langle l_{\sigma},x\rangle$ for $x\in\sigma$,
and $\langle l_\sigma, x\rangle=\langle l_\tau ,x\rangle$ whenever $x\in\tau<\sigma$.
We denote by $\SF(\Delta,\R)$ the additive group of support functions on $\Delta$.
\end{defn}
We will always assume that $|\Delta|$ is a convex cone.
A support function $h\in\SF(\Delta,\R)$ is said to be \emph{convex} if $h(x+y)\geq h(x)+h(y)$ for any
$x,y\in|\Delta|$.  We have for $\sigma\in\Delta(n)$, $\langle l_\sigma, x\rangle \geq h(x)$ for all $x\in|\Delta|$.
If for every $\sigma\in\Delta(n)$, we have equality only for $x\in\sigma$, then $h$ is said to be
\emph{strictly convex}.

Suppose $h\in\SF(\tilde{\Delta})$ is a strictly convex.  We will associate a rational convex polyhedral set
$\mathcal{C}_h \subset\mathfrak{t}^*$ to $\tilde{\Delta}$ and $h$.  Furthermore the fan associated to
$\mathcal{C}_h$ as in (\ref{eq:dual-cone-face}) is $\tilde{\Delta}$.
For each $\tau_j \in\tilde{\Delta}(1)$ we have a primitive element $u_j\in\Z_T,j=1,\ldots,N$, as above.
Set $\lambda_i :=h(u_i)$.  Then we define
\begin{equation}\label{eq:polyt-def}
\mathcal{C}_h :=\bigcap_{j=1}^{N} \{y\in\mathfrak{t}^* : \langle u_j,y\rangle\geq\lambda_j\}.
\end{equation}

We employ a construction originally due to Delzant and extended to the non-compact and singular cases by
D. Burns, V. Guillemin, and E. Lerman in~\cite{BurGuiLer} which constructs a K\"{a}hler structure on
$X_{\tilde{\Delta}}$ associated to a convex polyhedral set (\ref{eq:polyt-def}).  See also
\cite{Gui1,Gui2} for more on what is summarized here.
Let $\mathcal{A}:\Z^N \rightarrow\Z_T$ be the $\Z$-linear map with $\mathcal{A}(e_i)=u_i$, where $e_i,i=1,\ldots,N$,
are the standard basis vectors of $\Z^N$.  Then the $\R$-linear extension, also denoted by $\mathcal{A}$, induces a map
of Lie algebras $\mathcal{A}:\R^N \rightarrow\mathfrak{t}$.  Let $\mathfrak{k}=\ker\mathcal{A}$.  We have an exact
sequence
\begin{equation}\label{eq:lie-alg}
0\rightarrow\mathfrak{k}\overset{\mathcal{B}}{\longrightarrow}\R^N \overset{\mathcal{A}}{\longrightarrow}\mathfrak{t}\rightarrow 0,
\end{equation}
and its adjoint
\begin{equation}\label{eq:lie-alg-dual}
0\rightarrow\mathfrak{t}^* \overset{\mathcal{A}^*}{\longrightarrow}(\R^N)^* \overset{\mathcal{B}^*}{\longrightarrow}\mathfrak{k}^*\rightarrow 0.
\end{equation}
Also $\mathcal{A}$ induces a surjective map of Lie groups
$\ol{\mathcal{A}}:T^N \rightarrow T^n$, where $T^N =\R^{N}/{2\pi\Z^{N}}$.  If
$K=\ker\ol{\mathcal{A}}$, then we have the exact sequence
\begin{equation}\label{eq:lie-group}
1\rightarrow K\longrightarrow T^N \overset{\ol{\mathcal{A}}}{\longrightarrow} T^n \rightarrow 1.
\end{equation}
The moment map $\Phi$ for the action of $T^N$ on $(\C^N ,\frac{\mathbf{i}}{2}\sum_{j=1}^N dz_j \wedge d\ol{z}_j )$
is
\begin{equation}\label{eq:moment}
\Phi(z) =\sum_{j=1}^N |z_j |^2 e_j^*.
\end{equation}
Then moment map $\Phi_K$ for the action of $K$ on $\C^N$ is the composition
\begin{equation}
\Phi_K =\mathcal{B}^* \circ\Phi.
\end{equation}
Let $\lambda =\sum_{j=1}^N \lambda_j e_j^*$, and $\nu =\mathcal{B}^* (-\lambda)$.  Then
\begin{equation}
M_{\mathcal{C}_h} :=\Phi_K^{-1}(\nu)/K
\end{equation}
is smooth provided $\mathcal{C}_h$ in non-singular as in (\ref{eq:nonsing}).  The K\"{a}hler form on $\C^N$ descends
to a K\"{a}hler form $\omega_h$ on $M_{\mathcal{C}_h}$.  The action of
$T^n =T^N/K$ on $M_{\mathcal{C}_h}$ is Hamiltonian, and the restriction $\Phi|_{\Phi_K^{-1}(\nu)}$ descends to
$\ol{\Phi}:M_{\mathcal{C}_h}=\Phi_K^{-1}(\nu)/K \rightarrow(\R^N)^*$.  One can check that
$\im(\ol{\Phi}+\lambda)\subset\im(\mathcal{A}^*)$.  Thus
\begin{equation}
\Phi_{\mathcal{C}_h} :=(\mathcal{A}^*)^{-1}\circ (\ol{\Phi}+\lambda)
\end{equation}
is the moment map $\Phi_{\mathcal{C}_h}:M_{\mathcal{C}_h}\rightarrow \mathfrak{t}^*$ for $T^n =T^N/K$ acting on
$M_{\mathcal{C}_h}$.  Furthermore, $\im(\Phi_{\mathcal{C}_h}) =\mathcal{C}_h \subset\mathfrak{t}^*$.
Also, the action of $T^n$ on $M_{\mathcal{C}_h}$ extends to $T^n_{\C}\cong(\C^*)^n$.  And this action of
$T^n_{\C}$ has an open dense orbit.  Thus $M_{\mathcal{C}_h}$ is a toric variety.  And as the stablizer subgroups of
$T_{\C}^n$ coincide, $M_{\mathcal{C}_h}\cong X_{\tilde{\Delta}}$ as toric varieties.

We will make use of Guillemin's formula for the K\"{a}hler potential of $\omega_h$ on $M_{\mathcal{C}_h}$.
Let $l_j(y)=\langle u_j,y\rangle -\lambda_j$ for $j=1,\ldots,N$, and let
$l_{\infty}(y)=\sum_{j=1}^N \langle u_j,y\rangle$.  The following is proved in~\cite{BurGuiLer}; see also~\cite{Gui1,Gui2}.
\begin{thm}\label{thm:kahler-pot}
The K\"{a}hler form $\omega_h$ on the preimage $\Phi_{\mathcal{C}_h}^{-1}(\inter\mathcal{C}_h)$ of the
interior $\inter\mathcal{C}_h$ of the polyhedral set $\mathcal{C}_h$ is
\[\omega_h =\mathbf{i}\partial\ol{\partial}\Phi_{\mathcal{C}_h}^*(\sum_{j=1}^N \lambda_j \log(l_j) +l_{\infty}).\]
\end{thm}

Suppose $\tilde{\Delta}$ is a nonsingular subdivision of $\Delta$ giving a crepant resolution (\ref{eq:toric-resol}).
Then $u_1,\ldots,u_d \in\Z_T$ are vectors spanning the cone $\mathcal{C}^*(\mu)$, whereas
$u_{d+1},\ldots, u_N \in\Z_T$ are the lattice points in $\inter P_{\Delta}$.
In order to construct K\"{a}hler forms $\omega$ on $X_{\tilde{\Delta}}$ with $[\omega]\in H^2_c(X_{\tilde{\Delta}},\R)$
we make the following definition.
\begin{defn}
A strictly convex support function $h\in\SF(\tilde{\Delta},\R)$ is \emph{compact} if $h(u_j)=0$ for
$j=1,\ldots,d$.
\end{defn}

We will now prove Corollary~\ref{cor:main}.
If $h\in\SF(\tilde{\Delta},\R)$ is a compact strictly convex support function, then the K\"{a}hler form
$\omega_h$ on $X_{\tilde{\Delta}}$ has a compact K\"{a}hler class $[\omega_h ]$.
From Theorem~\ref{thm:kahler-pot} we have
\begin{equation}\label{eq:kahler-pot-com}
\omega_h =\mathbf{i}\partial\ol{\partial}\Phi_{\mathcal{C}_h}^*(\sum_{j=d+1}^N \lambda_j \log(l_j) +l_{\infty}).
\end{equation}
The potential function $F=\Phi_{\mathcal{C}_h}^*(\sum_{j=d+1}^N \lambda_j \log(l_j) +l_{\infty})$ is smooth
away from the exceptional set $E=\pi^{-1}(o)$, thus $[\omega_h ]\in H_c^2(X_{\tilde{\Delta}},\R)$.
We will construct a K\"{a}hler metric $\omega_0$ on $X_{\tilde{\Delta}}$ with all the properties in Lemma~\ref{lem:app-form}
which is furthermore invariant under $T^n$.  Let $f=\frac{r^2}{2}$ be the K\"{a}hler potential of the Ricci-flat
K\"{a}hler cone metric that exists by Theorem~\ref{thm:FOW}.  We consider $f$ as a function on
$X_{\tilde{\Delta}}$ via $\pi:X_{\tilde{\Delta}}\rightarrow X_{\Delta}$.  Let $0<a_1 <a_2$ and define a function
$\nu:\R_{>0}\rightarrow\R$ as in the proof of Lemma~\ref{lem:app-form}.  Then we have a non-negative form
$i\partial\ol{\partial}(\nu\circ f)\geq 0$ on $X_{\tilde{\Delta}}$ with $i\partial\ol{\partial}(\nu\circ f)>0$ on
$X_{a_2}=\{z\in X_{\tilde{\Delta}}: r(z)>a_2\}$.  Choose $b>a_2$, and choose $c_1>0$ large enough that
\begin{equation}
\Phi_{\mathcal{C}_h}^{-1}(\{y\in\mathcal{C}_h : l_\infty (y)\geq c_1 \})\subset X_b.
\end{equation}
Let $\phi:\R\rightarrow [0,1]$ be a smooth function with $\phi(x)=1$ for $x<c_1$ and $\phi(x)=0$ for $x>c_2$, where
$c_2 >c_1$.  Define $g=\phi\circ l_{\infty} \circ\Phi_{\mathcal{C}_h}$.  Then define
\begin{equation}\label{eq:app-form}
\omega_0 = i\partial\ol{\partial}(gF) +Ci\partial\ol{\partial}(\nu\circ f), \text{  for }C>0.
\end{equation}
For $C>0$ sufficiently large gives the metric with the required properties.

Corollary~\ref{cor:main} now follows from the proof of Theorem~\ref{thm:main}.  Since $\omega_0$ and the
Ricci-potential $f$ defined in (\ref{eq:ricci-pot}) are $T^n$-invariant, the $T^n$-invariance of the solution to
(\ref{eq:monge-amp}) follows from the uniqueness of the solution given in Proposition~\ref{prop:monge}.

Note that in Corollary~\ref{cor:main} we have a family of Ricci-flat K\"{a}hler metrics on $X_{\tilde{\Delta}}$
whose dimension is the number of lattice points in $\inter P_{\Delta}$, $N-d$ in the above notation.
For each $j=d+1,\ldots,N$, the prime divisor $D_j$ in $E=\pi^{-1}(o)$ is the smooth submanifold given by
$l_j \circ\Phi_{\mathcal{C}_h} =0$.  Let $c_j \in H^2_c (X_{\tilde{\Delta}},\R)$ be the cohomology dual of
$[D_j]$ in $H_{2n-2}(X_{\tilde{\Delta}},\R)$.   Then $c_j =[\beta_j]$ (cf.~\cite{Gui1}, Theorem 6.2), where
\begin{equation}
\beta_j =\frac{i}{2\pi}\ol{\partial}\partial\log\Phi_{\mathcal{C}_h}^* l_j.
\end{equation}
If $\omega$ is the K\"{a}hler form of Corollary~\ref{cor:main} starting with $\omega_0$ in (\ref{eq:app-form}),
then we have from (\ref{eq:kahler-pot-com}) that
\begin{equation}
[\omega] =-2\pi\sum_{j=d+1}^N \lambda_j c_j.
\end{equation}

\section{Examples}

\subsection{Asymptotically locally Euclidean K\"{a}hler manifolds}

Let $\Gamma\subset GL(n,\C)$ be a finite subgroup, and consider the singular space $X=\C^n /\Gamma$.
We want isolated singularities, so we assume $\Gamma$ acts freely on $\C^n \setminus\{o\}$.
The singularity is Gorenstein precisely when $\Gamma\subset SL(n,\C)$.  We may assume $\Gamma\subset SU(n)$,
as $\Gamma$ is always conjugate to such a subgroup.  Note that this is precisely the case in which
$S=S^{2n-1} /\Gamma$ has constant curvature and $C(S)=X\setminus\{o\}$ is flat.

When $n=2$, $X=\C^2 /\Gamma$ is a \emph{Kleinian singularity}.  And $X$ admits a unique crepant resolution
$\pi: Y\rightarrow X$.  For $n=3$, $X=\C^3 /\Gamma$, it was proved by S. Roan~\cite{Ro} that $X$ admits a crepant
resolution, but it may not be unique.  For $n\geq 4$, $X$ may or may not admit a crepant resolution, and if it exists
it may or may not be unique.

In this case $H_c^2(Y,\R)=H^2(Y,\R)$, Theorem~\ref{thm:main} shows that there is a Ricci-flat K\"{a}hler metric in
every K\"{a}hler class asymptotic to the flat metric as in (\ref{eq:asymp}).  But in this case there is an improved proof,
by D. Joyce~\cite{Joy1,Joy2}, which proves Theorem~\ref{thm:main} where one has (\ref{eq:asymp}) with $\delta=0$.

\subsection{Canonical bundles of toric Fano manifolds}

Let $M$ be a Fano manifold.  Then the canonical bundle $\mathbf{K}_M$ is negative.   Let $Y=\mathbf{K}_M$ denote
the total space.  Then we have the Remmert reduction (cf.~\cite{Gra}) $\pi: Y\rightarrow X$ which collapses the
zero section of $\mathbf{K}_M$.  Then $X =C(S)\cup\{o\}$, where $S$ the $U(1)$-subbundle of $\mathbf{K}_M$ with
the usual Sasaki structure.  It is not difficult to check that $\pi: Y\rightarrow X$ is a crepant resolution.
If $M$ admits a K\"{a}hler-Einstein metric, then after a possible $D$-homothetic transformation as in (\ref{eq:d-homo}),
the standard Sasaki structure on $S$ as in Example~\ref{xpl:Sasak-st} is Sasaki-Einstein.
The Calabi ansatz gives a complete Ricci-flat K\"{a}hler
metric on $Y$ (cf.~\cite{Cal}).  If $M$ is not K\"{a}hler-Einstein, then $S$ can possibly have a Sasaki-Einstein
structure for a different Reeb vector field.

Suppose $M$ is a toric Fano manifold of dimension $m$.  We have a crepant resolution $\pi:Y\rightarrow X$ as above,
where $X=C(S)$ is a toric K\"{a}hler cone satisfying Proposition~\ref{prop:CY-cond-toric}.
By Theorem~\ref{thm:FOW} $X=C(S)$ has a Ricci-flat K\"{a}hler cone metric for some Reeb vector field.  And
by Corollary~\ref{cor:main} there is a 1-dimensional family of asymptotically conical Ricci-flat K\"{a}her metrics on
$Y$.

We have that $M$ is given by a fan $\Delta$ in $\Z^m$, and we give the fans $\ol{\Delta}$ of $X=C(S)$ and
$\tilde{\Delta}$ of $Y=\mathbf{K}_M$ in $\Z^n$, $n=m+1$.
If $u_1,\ldots,u_d \in\Z^n$ are primitive elements generating each $\tau\in\Delta(1)$, then $\ol{\Delta}$ consists of
the convex polyhedral cone spanned by $\ol{u}_1 =(u_1, 1),\ol{u}_2 =(u_2 ,1),\ldots,\ol{u}_d =(u_d ,1)$ and all of
its faces.

Let $\alpha$ be the 1-cone generated by $e_n \in\Z^n$.  Then $\tilde{\Delta}$ consists of all cones of $\ol{\Delta}$
besides the n-dimensional cone plus the following.  For $\sigma\in\ol{\Delta}(r), r<n$, let
$\tilde{\sigma} =\sigma + \alpha$.  It is easy to see that this defines a non-singular subdivision of $\ol{\Delta}$.

\begin{figure}[htb]
 \centering
 \includegraphics[scale=0.4]{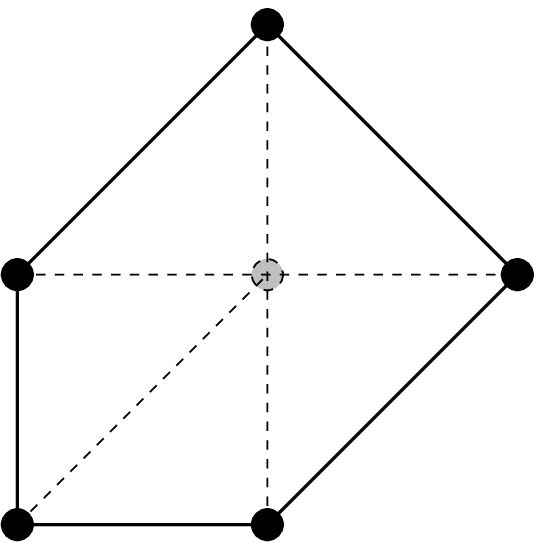}
 \caption{Canonical bundle of $\cps^2_{(2)}$}
 \label{fig1}
\end{figure}

Consider $M=\cps_{(2)}^2$, the two-points blow up.  Then $X=C(S)$ has a Ricci-flat K\"{a}hler cone metric as in
Example~\ref{xpl:two-points}, for a non-regular Sasaki-Einstein structure.  The lattice triangulation of
the polytope $P_{\ol{\Delta}}$ is given in Figure~\ref{fig1}.

\subsection{Resolutions of $C(Y^{p,q})$}

A series of 5-dimensional Sasaki-Einstein metrics $Y^{p,q}$, with $p,q\in\N, p>q>0$, and $\gcd(p,q)=1$, first appeared
in~\cite{GMSW1}.  These examples are remarkable in that they contain the first known examples of irregular
Sasaki-Einstein manifolds, and also because the metrics are given explicitly.  These examples are toric and are further of
cohomogeneity one with an isometry group of $SO(3)\times U(1)\times U(1)$ if $p,q$ are both odd, and
$U(2)\times U(1)$ otherwise.

The Sasaki structure is quasi-regular precisely when $p,q\in\N$ as above satisfy the diophantine equation
\begin{equation}
4p^2 -3q^2 =r^2,
\end{equation}
for some $r\in\Z$.  It was shown in~\cite{GMSW1} that there are both infinitely many quasi-regular and irregular
examples.

\begin{figure}[htb]
 \centering
 \includegraphics[scale=0.4]{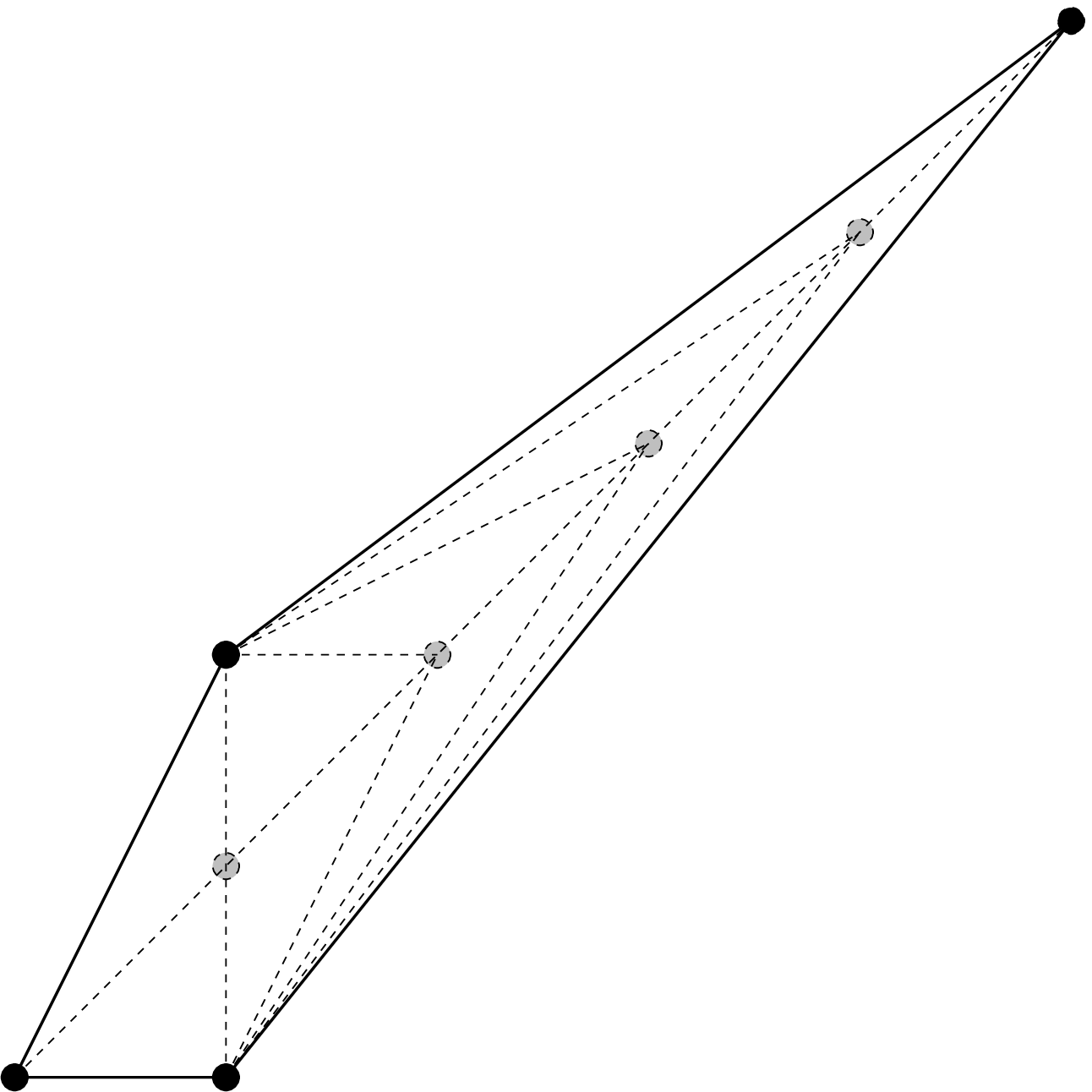}
 \caption{A resolution of $X^{5,3}$}
 \label{fig2}
\end{figure}

We have $X_{\Delta} =C(Y^{p,q})\cup\{o\}$ where the fan $\Delta$ in $\Z^3$ is generated by the four vectors
\begin{equation}
u_1 =(0,0,1), u_2 =(1,0,1), u_3 =(p,p,1), u_4 =(p-q-1,p-q,1).
\end{equation}
A basic lattice triangulation of $P_{\Delta}$ can be constructed for general $p,q$ as is shown in Figure~\ref{fig2} for
$Y^{5,3}$.  It is not difficult to see that the subdivision $\tilde{\Delta}$ of $\Delta$ has a compact strictly
convex support function.  Thus Corollary~\ref{cor:main} gives a $p-1$-dimensional family of asymptotically conical
Ricci-flat K\"{a}hler metrics on $X_{\tilde{\Delta}}$.

\subsection{Toric crepant resolutions}

\begin{figure}[htb]
 \centering
 \includegraphics[scale=0.4]{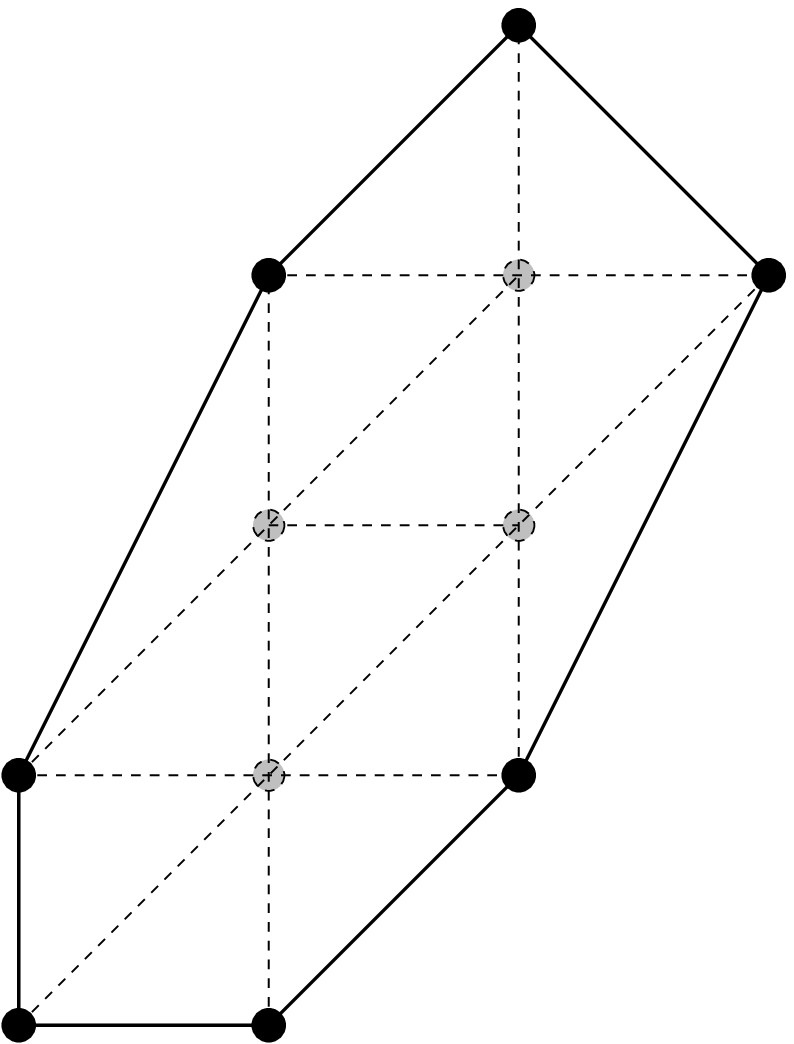}
 \caption{Example}
 \label{fig3}
\end{figure}

\begin{figure}[t]
 \centering
 \includegraphics[scale=0.4]{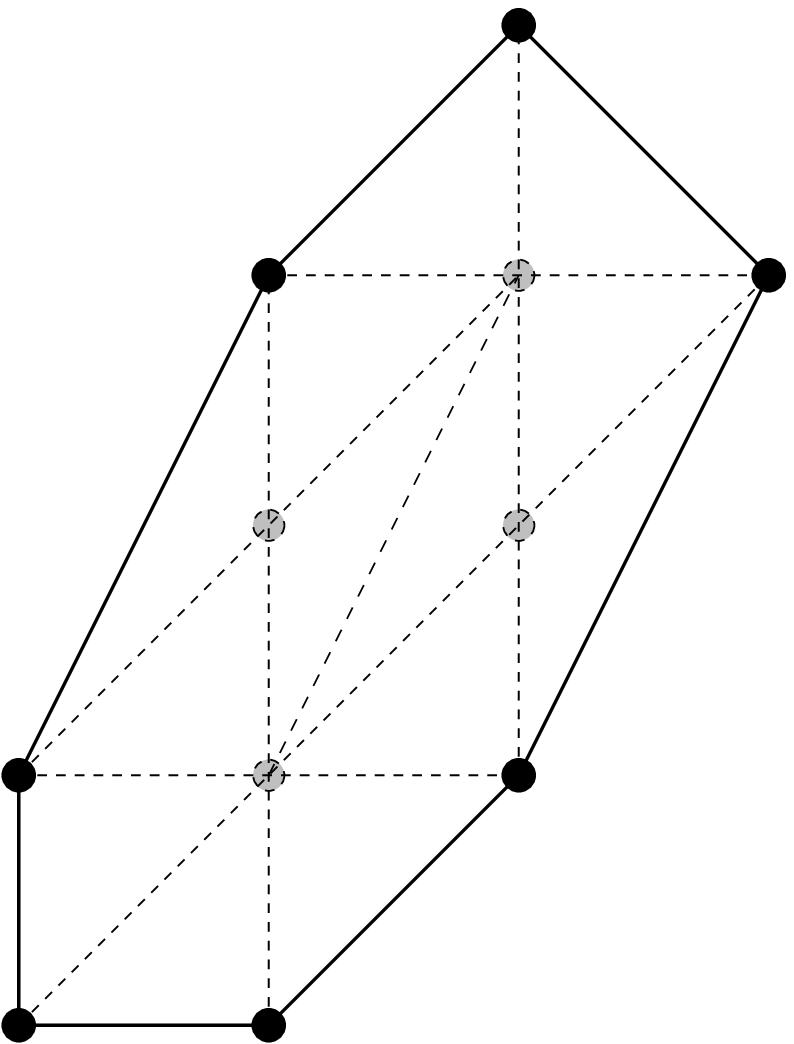}
 \caption{a \emph{flop} of~\ref{fig3}}
 \label{fig4}
\end{figure}

Let $X_{\Delta}$ be a toric K\"{a}hler cone.  If $n=3$ then Proposition~\ref{prop:3dim-resol} implies that $X_{\Delta}$
admits a toric crepant resolution, $\pi:X_{\tilde{\Delta}}\rightarrow X_{\Delta}$.  And more generally, if $n>3$, then
$X_{\Delta}$ admits a toric partial crepant resolution $\pi:X_{\tilde{\Delta}}\rightarrow X_{\Delta}$ which has at most orbifold singularities.  The author does not have a general result on the existence of a compact strictly convex support function on
$\tilde{\Delta}$.  Nevertheless, it is elementary to construct examples, such as in Figure~\ref{fig3}, which has a
4-dimensional space of asymptotically conical Ricci-flat K\"{a}hler metrics.  In this example $X_{\Delta}$ has
another resolution, Figure~\ref{fig4}, which is related to Figure~\ref{fig3} by a flop.

It is proved in~\cite{vC3} that for $n=3$, as long as $X$ is not the quadric cone,
there is a crepant resolution $\pi:X_{\tilde{\Delta}}\rightarrow X_{\Delta}$ such that $\tilde{\Delta}$ has a compact
strictly upper convex support function.  And therefore, Corollary~\ref{cor:main} applies.
This can be used to easily construct infinitely many 3-dimensional examples.

\pagebreak

\bibliographystyle{plain}

\end{document}